\documentclass{scrartcl}
\usepackage{amsmath}
\usepackage{amsfonts}
\usepackage{amssymb}
\usepackage{dsfont}
\usepackage[T1]{fontenc}
\usepackage[latin1]{inputenc}
\usepackage{epsfig}
\usepackage{graphicx}

\usepackage[boxed, lined]{algorithm2e}
\usepackage{float}
\usepackage{comment}
\usepackage{mathtools}

\usepackage{bm}

\newcommand{\subalign}[1]{%
   \vcenter{%
     \Let@ \restore@math@cr \default@tag
     \baselineskip\fontdimen10 \scriptfont\tw@
     \advance\baselineskip\fontdimen12 \scriptfont\tw@
     \lineskip\thr@@\fontdimen8 \scriptfont\thr@@
     \lineskiplimit\lineskip
     \ialign{\hfil$\m@th\scriptstyle##$&$\m@th\scriptstyle{}##$\hfil\crcr
       #1\crcr
     }%
   }%
}

\newcommand{\bk}{\mathbf{k}}

\newcommand{\bl}{\mathbf{l}}
\newcommand{\ba}{\mathbf{a}}
\newcommand{\bb}{\mathbf{b}}

\newcommand{\be}{\mathbf{e}}
\newcommand{\br}{\mathbf{r}}

\newcommand{\bw}{\mathbf{w}}
\newcommand{\bv}{\mathbf{v}}

\newcommand{\D}{{\mathcal{D}}}

\newcommand{\W}{{\mathcal{W}}}

\newcommand{\Nu}{{\mathcal{N}}}

\newcommand{\N}{\mathbb{N}}

\newcommand{\R}{\mathbb{R}}

\newcommand{\Rd}{\mathbb{R}^d}

\newcommand{\beq}{\begin{eqnarray*}}

\newcommand{\eeq}{\end{eqnarray*}}

\newcommand{\beqm}{\begin{eqnarray}}

\newcommand{\eeqm}{\end{eqnarray}}

\newtheorem{theorem}{Theorem}

\newtheorem{lemma}{Lemma}
\newtheorem{definition}{Definition}

\DeclareOldFontCommand{\bf}{\normalfont\bfseries}{\mathbf}
\DeclareOldFontCommand{\it}{\normalfont\itshape}{\mathit}

\newcommand{\EXP}{{\mathbf E}}
\newcommand{\PROB}{{\mathbf P}}

\renewcommand{\P}{{\cal P}}
\renewcommand{\bf}{\normalfont \bfseries}
\renewcommand{\it}{\normalfont \itshape}

\allowdisplaybreaks
\begin{document}
\renewcommand{\thefootnote}{\fnsymbol{footnote}}
\newcommand{\F}{{\cal F}}
\newcommand{\Sp}{{\cal S}}
\newcommand{\G}{{\cal G}}
\newcommand{\HH}{{\cal H}}

\begin{center}

  {\LARGE \bf
    On the rate of convergence of an over-parametrized deep
    neural network regression estimate learned by gradient
    descent
  }
\footnote{
Running title: {\it Neural network estimate learned by gradient descent}}
\vspace{0.5cm}

Michael Kohler

{\it 
Fachbereich Mathematik, Technische Universit\"at Darmstadt,
Schlossgartenstr. 7, 64289 Darmstadt, Germany,
email: kohler@mathematik.tu-darmstadt.de.
}

\end{center}
\vspace{0.5cm}

\begin{center}
April 4, 2025
\end{center}
\vspace{0.5cm}

\noindent
    {\bf Abstract}\\

    Nonparametric regression with random design is considered.
    The $L_2$ error with integration with respect to the design
    measure is used as the error criterion.
    An over-parametrized deep neural network
    regression estimate
    with logistic activation function
    is defined, where all weights are learned
    by gradient descent. It is shown that the estimate
    achieves a nearly optimal rate of convergence in case
    that the regression function is $(p,C)$--smooth.
    
    \vspace*{0.2cm}


\vspace*{0.2cm}


\section{Introduction}
\label{se1}
\subsection{Scope of this paper}
Deep learning has achieved tremendous success in applications,
e.g., in
image classification
(cf., e.g., Krizhevsky, Sutskever and Hinton  (2012)),
language recognition (cf., e.g.,  Kim (2014))
machine translation (cf., e.g., Wu et al. (2016)),
mastering of games (cf., e.g., Silver et al. (2017))
or simulation of human conversation (cf., e.g., Zong and Krishnamachari
(2022)).
From a theoretical point of view this great success is still
a mystery. In particular, it is unclear why
the use of over-parametrized
deep neural networks, which have much more weights than there are data
points, does not lead to an overfitting of the estimate, and why
gradient descent is able
to minimize the nonlinear and non-convex empirical risk in the
definition of the estimates in such a way that the estimates can
achieve a small risk for randomly initialized starting values
for the weights of the networks. In a standard
regression setting we give answers to these two questions
in the special situation that we want to estimate a
$(p,C)$--smooth regression function.

\subsection{Nonparametric regression}
We study deep neural networks in the context of nonparametric
regression. Here we have given an $\R^d \times \R$--valued random
vector $(X,Y)$ with $\EXP Y^2 < \infty$, and our goal is to
predict the value of $Y$ given the value of $X$. Let $m(x)=\EXP\{Y|X=x\}$
be the so--called regression function. Then any measurable $f:\R^d \rightarrow \R$ satisfies
\begin{equation}
  \label{se1eq1}
\EXP \{|f(X)-Y|^2\} = \EXP\{|m(X)-Y|^2\} + \int |f(x)-m(x)|^2 \PROB_X (dx)
\end{equation}
(cf., e.g., Section 1.1 in Gy\"orfi et al. (2002)), hence in view of
minimizing the so-called $L_2$ risk (\ref{se1eq1}) of $f$ the regression function $m$ is the optimal predictor, and the so--called $L_2$ error
\begin{equation}
  \label{se1eq2}
\int |f(x)-m(x)|^2 \PROB_X (dx)
\end{equation}
describes how far the $L_2$ risk of a function $f$ is away from its optimal value.

In applications typically the distribution of $(X,Y)$ and hence also
the corresponding regression function $m$ is unknown. But often it is
possible to observe data from the underlying distribution, and the task
is to use this data to estimate the unknown regression function.
In view of minimization of the $L_2$ risk of the estimate, here it is
natural to use the $L_2$ error as an error criterion.

In order to introduce this problem formally, let $(X,Y)$, $(X_1,Y_1)$,
\dots, $(X_n,Y_n)$ be independent and identically distributed. In nonparametric
regression the data set
\begin{equation}
  \label{se1eq3}
  \D_n = \left\{
(X_1,Y_1), \dots, (X_n,Y_n)
  \right\}
\end{equation}
is given, and the task is to construct an estimate
\[
m_n(\cdot)=m_n(\cdot,\D_n):\R^d \rightarrow \R
\]
such that its $L_2$ error
\[
\int | m_n(x)-m(x)|^2 \PROB_X (dx)
\]
is small. A systematic introduction to nonparametric regression,
its estimates and known results can be found, e.g., in
Gy\"orfi et al. (2002).

\subsection{Rate of convergence}
Stone (1982) determined the optimal Minimax rate of convergence
of the $L_2$ error in case of a smooth regression function.
Here he considered so-called $(p,C)$--smooth regression functions,
which are defined as follows.

\begin{definition}
\label{se1de1} 
  Let $p=q+s$ for some $q \in \N_0$ and $0< s \leq 1$.
A function $m:\R^d \rightarrow \R$ is called
$(p,C)$-smooth, if for every $\bm{\alpha}=(\alpha_1, \dots, \alpha_d) \in
\N_0^d$
with $\sum_{j=1}^d \alpha_j = q$ the partial derivative
$\partial^q m/(\partial x_1^{\alpha_1}
\dots
\partial x_d^{\alpha_d}
)$
exists and satisfies
\[
\left|
\frac{
\partial^q m
}{
\partial x_1^{\alpha_1}
\dots
\partial x_d^{\alpha_d}
}
(x)
-
\frac{
\partial^q m
}{
\partial x_1^{\alpha_1}
\dots
\partial x_d^{\alpha_d}
}
(z)
\right|
\leq
C
\cdot
\|\bold{x}-\bold{z}\|^s
\]
for all $\bold{x},\bold{z} \in \R^d$, where $\Vert\cdot\Vert$ denotes the Euclidean norm.
\end{definition}

Stone (1982) showed that in case of a $(p,C)$--smooth regression function the
optimal Minimax rate of convergence for the expected $L_2$ error is
\[
n^{-\frac{2p}{2p+d}}.
\]
This rate suffers from the so--called curse of dimensionality: If
the dimension $d$ is large
compared to the smoothness $p$ of the regression function, the exponent
will be close to zero and the rate of convergence will be rather
slow. Since this rate is optimal, the only way to circumvent this
is to impose additional assumptions on the structure of the regression function.
For this, various assumptions exists, e.g., additive models (cf., e.g.,
Stone (1985)),
interaction models (cf., e.g., Stone (1994)),
single index models (cf., e.g.,
 H\"ardle, Hall and Ichimura (1993), H\"ardle and Stoker (1989),
Yu and Ruppert (2002) and Kong and Xia (2007))
or projection pursuit (cf, e.g.,
Friedman and Stuetzle (1981)), where corresponding low dimensional
rates of convergence can be achieved (cf., e.g., Stone (1985, 1994)
and Chapter 22 in Gy\"orfi et al. (2002)).

\subsection{Least squares neural network estimates}
One way to estimate a regression function is to define a function
space $\F_n$ consisting of functions $f: \R^d \rightarrow \R$ and
to use the principle of least squares to select one of its functions as
the regression estimate, i.e., to define
\[
m_n(\cdot) = \arg \min_{f \in \F_n}
\frac{1}{n} \sum_{i=1}^n |f(X_i)-Y_i|^2.
\]
In view of the $L_2$ error of the estimate it is important
that the function space is on the one hand large enough such that
a  function is contained in it which approximates the unknown
regression function well, and that on the other hand the
function space is not too complex so that the empirical $L_2$ risk
\[
\frac{1}{n} \sum_{i=1}^n |f(X_i)-Y_i|^2
\]
is not too far from the $L_2$ risk for functions in this function space.

One possible way to define such function spaces in case of $d$ large is
to use feedforward neural networks. These function depend on an
activation function $\sigma: \R \rightarrow \R$, e.g.,
\[
\sigma(x)=\max \{x,0\}
\]
(so-called ReLU-activation function) or
\[
\sigma(x)= \frac{1}{1+e^{-x}}
\]
(so-called logistic activation function).

The most simple form of neural networks are shallow networks, i.e., neural networks
with one hidden layer, in which a simple linear combination
of artifical neurons defined by applying the activation function
to a linear combination of the components of the input is used
to define a function $f:\Rd \rightarrow \R$
by
\begin{equation}
  \label{se1eq3b}
f(x)= \sum_{k=1}^K \alpha_k \cdot
\sigma \left(
\sum_{j=1}^d \beta_{k,j} \cdot x^{(j)} + \beta_{k,0}
\right)
+ \alpha_0.
\end{equation}
Here $K \in \N$ is the number of neurons, and the weights
$\alpha_k \in \R$ $(k=0, \dots, K)$, $\beta_{k,j} \in \R$
$(k=1, \dots, K, j=0, \dots, d)$ are chosen by the principle of
least squares.

The rate of convergence
of shallow neural networks regression estimates
has been analyzed in
Barron (1994) and McCaffrey and Gallant (1994).
Barron (1994) proved a dimensionless rate of $n^{-1/2}$
(up to some logarithmic factor), provided the Fourier transform
of the regression function
has a finite first
moment, which basically
requires that the function becomes smoother with increasing
dimension $d$ of $X$.
McCaffrey and Gallant (1994) showed a rate of $n^{-\frac{2p}{2p+d+5}+\varepsilon}$
in case of a $(p,C)$-smooth regression function, but their study was restricted to the use of a certain cosine squasher as activation function.

In deep learning neural networks with
several hidden layers are used to define classes of functions.
Here
a feedforward neural network with $L \in \N$ hidden layers and $k_s \in \N$
neurons in the layers $s \in \{1, \dots, L\}$ is recursively defined
by
\begin{equation}
  \label{se1eq4}
  f_\bw(x)
=
\sum_{j\in \{1, \dots,k_L\}}
w_{1,j}^{(L)} \cdot f_j^{(L)}(x),
\end{equation}
where
\begin{equation}
\label{se1eq5}
f_i^{(s)}(x)
=
\sigma \left(
\sum_{
j \in \{1, \dots,k_{s-1}\}}
w_{i,j}^{(s-1)} \cdot f_j^{(s-1)}(x)
+ w_{i,0}^{(s-1)}
\right)
\quad \mbox{for }
s \in \{2, \dots, L\} \mbox{ and } i>0
\end{equation}
and
\begin{equation}
\label{se1eq6}
f_i^{(1)}(x)
=
\sigma \left(
\sum_{
j \in \{1, \dots, d\}
}
w_{i,j}^{(0)} \cdot x^{(j)} + w_{i,0}^{(0)} 
\right)
\quad \mbox{for }
i>0.
\end{equation}

The rate of convergence of least squares estimates
based on multilayer neural networks 
has been analyzed in Kohler and Krzy\.zak (2017),
 Imaizumi and Fukamizu (2018),
Bauer and Kohler (2019),
Suzuki and Nitanda (2019),
Schmidt-Hieber (2020) and Kohler and Langer (2021).
One of the main results achieved in this context shows
that neural networks
can achieve some kind of dimension reduction under
rather general assumptions. The most general form goes back to
Schmidt-Hieber (2020) and can be formalized as follows:

\begin{definition}
\label{se1de2}
Let $d \in \N$ and $m: \Rd \to \R$ and let
$\P$ be a subset
of $(0,\infty) \times \N$.\\
\noindent
\textbf{a)}
We say that $m$ satisfies a hierarchical composition model of level $0$
with order and smoothness constraint $\mathcal{P}$, if there exists a $K \in \{1, \dots, d\}$ such that
\[
m(\bold{x}) = x^{(K)} \quad \mbox{for all } \bold{x} = (x^{(1)}, \dots, x^{(d)})^{\top} \in \Rd.
\]
\noindent
\textbf{b)}
We say that $m$ satisfies a hierarchical composition model
of level $l+1$ with order and smoothness constraint $\mathcal{P}$, if there exist $(p,K)  \in \P$, $C>0$,
$g: \R^{K} \to \R$ and $f_{1}, \dots, f_{K}: \Rd \to \R$, such that
$g$ is $(p,C)$--smooth,
$f_{1}, \dots, f_{K}$ satisfy a  hierarchical composition model of level $l$
with order and smoothness constraint $\mathcal{P}$
and 
\[m(\bold{x})=g(f_{1}(\bold{x}), \dots, f_{K}(\bold{x})) \quad \mbox{for all } \bold{x} \in \Rd.\]
\end{definition}

Schmid-Hieber (2020) showed that suitable least squares neural network 
regression estimates
achieve (up to some logarithmic factor) a rate of convergence of order
\[
\max_{(p,K) \in \mathcal{P}} n^{- \frac{2p}{2p+K}}
\]
in case that the regression function satisfies a hierarchical composition
model of some finite level with order and smoothness constraint $\mathcal{P}$.
Since this rate of convergence does not depend on the dimension $d$ of $X$,
this results shows that least squares neural network regression estimates
are able to circumvent the curse of dimensionality in case that the
regression function satisfies a hierarchical composition model.

\subsection{Learning of neural network estimates}
In applications the least squares estimates of the previous
subsection cannot be used, since it is not clear how one can
minimize the empirical $L_2$ risk, which is a nonlinear and
non-convex function of the weights. Instead, one uses
gradient descent applied to a randomly chosen starting vector of weights
to minimize it approximately. Typically, here the estimates are
over-parameterized, i.e., they use much more weights than there
are data points, so in principle it is possible to choose the weights
such that the data points are interpolated (at least, if the x values
are all distinct).

In practice it has been observed, that this procedure
leads to estimates which predict well on new independent test data.
There have been various attemps to explain this using some models
for deep learning. E.g., Choromanska et al. (2015) used random matrix
theory to show that in this model of deep learning
the so--called landscape hypothesis is true,
which states that the loss surface contains many deep local minima.
Other popular models for deep learning include the neural tangent
kernel setting proposed by Jacot, Gabriel and Hongler (2020)
or the meanfield approach (cf., e.g., Mei, Montanari, and Nguyen (2018)).
The problem with studying deep neural networks in equivalent models is
that it is unclear how close the behaviour of the deep networks in
the proposed
equivalent model is to the behaviour of the deep networks in the
applications, because they are based on some approximation of the
application using e.g. some asymptotic expansions.

There exits also various articles which study over-parametrized
deep neural network estimates directly in a standard regression
setting. Kohler and Krzy\.zak (2022) showed that these estimates
can achieve a nearly optimal rate of convergence in case that
the regression function is $(p,C)$--smooth with $p=1/2$. Furthermore
it was shown there that these estimates can be modified such that
they achieve a dimension reduction in an interaction model. These
results require that a penalized empirical $L_2$ risk is minimized
by gradient descent. That such results also hold without using
any regularization by a penalty term was shown in Drews and Kohler (2023).
Again, the estimates achieve a nearly optimal rate of convergence
only in case of a $(p,C)$--smooth regression function with $p=1/2$.

\subsection{Main results}
In this article we extend the results from Kohler and Krzy\.zak (2022)
and Drews and Kohler (2023) to the case of a general $p \in [1/2, \infty)$.
To do this, we
study over-parametrized
deep neural network regression estimates with logistic activation function,
where the values of all the weights
are learned by gradient descent, in a standard regression model. 
We consider a $(p,C)$--smooth regression function,
and we choose as topology of the neural network a linear combination of
many parallel computed deep neural networks of fixed depth and width.
We show that for a suitable initialization of the weights, a suitably
chosen
stepsize of the gradient descent, and a suitably chosen
number of gradient descent steps the expected $L_2$ error of
our estimates converges to zero with rate
\[
n^{- \frac{2p}{2p+d} + \epsilon},
\]
where $\epsilon>0$ can be chosen as an arbitrary small constant.

In order to prove this result we show three crucial auxiliary results:
Firstly,
we show that during gradient descent our estimates stay in a function space
which has a finite complexity (measured by its supremum norm covering number).
We achieve this by showing that the
weights remain bounded and consequently the derivatives of the estimate
stay bounded, which enables us to bound the covering number using metric
entropy bounds.
Secondly, we derive new approximation results for neural networks
with bounded weights, where the bounds fit the upper bounds
on the covering number derived by using the metric entropy bounds.
And thirdly, we show that the gradient descent is linked to a
gradient descent applied to a linear Taylor polynomial of our network,
and therefore can be analyzed by techniques develloped for
the analysis of gradient descent for smooth convex functions.

In our theory over-parametrized deep neural networks do not overfit
because the weights remain bounded during training and consequently the
networks stay in a function space of bounded complexity.
Furthermore, the gradient descent can find neural networks which approximate
the unknown regression function well since the over-parametrized
structure and the initialization
of our network is such that with high probability there
is a network with good approximation properties close to our
initial network.

\subsection{Discussion of related results}

Motivated by the huge success of deep learning in applications,
there have been already quite a few results derived concerning the
theoretical analysis of these methods. E.g., there exist many results
in approximation theory for deep neural networks, see, e.g.,
Yarotsky (2018),
Yarotsky and Zhevnerchute (2019),
Lu et al. (2020), Langer (2021) and the literature cited therein.
These results show that smooth functions can be approximated
well by deep neural networks and analyze what kind of topology
and how many nonzero weights are necessary to approximate
a smooth function up to some given error. In applications,
the functions which one wants to approximate has to be estimated
from observed data, which usually contains some random error.
It has been also already analyzed how well a network learned
from such noisy data generalizes on new independent test data.
This has been done within the framework of the classical VC theory
(using e.g. the result of
Bartlett et al. (2019) to bound the VC dimension of classes
of neural networks) or in case of over-parametrized deep
neural networks
(where the number of free parameters adjusted
to the observed data set is much larger than the sample size)
by using bounds on the Rademacher complexity
(cf., e.g., Liang, Rakhlin and Sridharan (2015), Golowich, Rakhlin and
Shamir (2019), Lin and Zhang (2019),
Wang and Ma (2022)
and the literature cited therein).
By combining these kind of results it was possible to analyze the
error of least squares regression estimates. Here it was shown in
a series of papers (cf., e.g., Kohler and Krzy\.zak (2017), Bauer and
Kohler (2019), Schmidt-Hieber (2020) and Kohler and Langer (2021))
that least squares regression estimates based on deep networks
can achieve a dimension reduction in case that the function to be
estimates satisfies a hierarchical composition model, i.e., in case
that it is a composition of smooth functions which do either depend
only on a few components or are rather smooth. This is due to the
network structure of deep networks, which implies that the composition
of networks is itself a deep network. Consequently, any
approximation
result of some kind of functions by deep networks can be extended 
to an approximation result of a composition of such function
by a deep network representing a composition of the approximating
networks. And hereby the number of weights and the depth of the
network, which determine the VC dimension and hence the
complexity of the network in case that it is not over-parametrized
(cf., Bartlett et al. (2019)), changes not much. So such a network
has the approximation properties and the complexity of a network
for low dimensional predictors and hence can achieve a dimension
reduction.

There also exist quite a few results on the optimization of deep
neural networks. E.g., Zou et al. (2018), Du et al. (2019),
Allen-Zhu, Li and Song (2019) and Kawaguchi and Huang (2019)
analyzed the
 application of gradient descent to over-parameterized deep
neural networks. It was shown in these papers that this
leads to neural networks which (globally) minimize the empirical risk
considered. Unfortunately, as was shown in  Kohler and Krzy\.zak (2021),
the corresponding estimates do not behave well on new independent
data.

As pointed out by  Kutyniok (2020), it is essential for
a theoretical analysis of deep learning estimates to study
simultaneously
 the approximation
 error, the generalization error and the optimization error,
 and none of the results mentioned above
controlls all these three aspects together.

There exists various approaches where these three things are studied
simultaneously in some equivalent models of deep learning.
The most prominent  approach here is the
neural tangent kernel setting proposed by Jacot, Gabriel and Hongler (2020).
Here instead of a neural network
estimate a kernel estimate is studied and its error is used
to bound the error of the neural network estimate. For further
results in this context see Hanin and Nica (2019) and the literature
cited therein.
As was pointed out in Nitanda and Suzuki (2021) in most studies in the
 neural tangent kernel setting the equivalence to deep neural networks
holds only pointwise and not for the global $L_2$ error, hence from
these
result it is not clear how the $L_2$ error of the deep neural network
estimate behaves.
Nitanda and Suzuki (2021) were able to analyze the global error
of 
an over-parametrized shallow neural network
learned by gradient descent based on this approach. However, due to the
use of the neural tangent kernel, also the smoothness assumption
of the function to be estimated has to be defined with the aid of
a norm involving the kernel, which does not lead to classical
smoothness conditions, which makes it hard to
understand the meaning of the results. 
Furthermore, their result did not specify how many neurons the
shallow neural network must have, it was only shown that the results
hold if this number of neurons is sufficiently large, and it is not
clear
whether it must grow, e.g., exponentially in the sample size or not.
Another approach where the estimate
is studied  in some asymptotically equivalent model
is the  mean field approach, cf.,
Mei, Montanari, and Nguyen (2018), Chizat and Bach (2018) or Nguyen and Pham (2020).

The theory presented in this article is an extension of the theory
develloped in Braun et al. (2023), Drews and Kohler (2023, 2024)
and Kohler and Krzy\.zak (2022, 2023). The basic idea there is that
for smooth activation functions the inner weights do not change
much during learning if the stepsizes are sufficiently small and
it was shown that at the same time the outer weights will be chosen
suitably by gradient descent. In this article we extend this theory
by showing that in our special topology gradient descent is also able
to learn the inner weights locally, and by deriving a new
approximation result for the approximation of $(p,C)$--smooth
functions by deep neural network with bounded weights. In fact,
it is the new approximation results which is essential to
extend the previous results from $(p,C)$--smooth regression functions
with $p=1/2$ to the case of  $(p,C)$--smooth regression function
with  general $p \geq 1/2$, and it can be shown that  the
rate of convergence of this article can also be achieved
if only the weights of the output
layer are changed during gradient descent and all other weights
keep their initial values (cf., Remark 1).
This approach is related to the so--called
random feature networks, where the inner weights are not learned
at all and gradient descent is applied only to the weights in the
output level, cf., e.g., Huang, Chen and Siew (2006) and
Rahimi and Recht (2008a, 2008b, 2009).

\subsection{Notation}
\label{se1sub5}
  The sets of natural numbers, real numbers and nonnegative real numbers
  are denoted by $\N$, $\R$ and $\R_+$, respectively.
  For $z \in \R$, we denote
the smallest integer greater than or equal to $z$ by
$\lceil z \rceil$.
The Euclidean norm of $x \in \Rd$
is denoted by $\|x\|$. For a closed and convex set $A \subseteq \R^d$
we denote by $Proj_A x$ that element $Proj_A x \in A$ with
\[
\|x-Proj_A x\|= \min_{z \in A} \|x-z\|.
\]
For $f:\R^d \rightarrow \R$
\[
\|f\|_\infty = \sup_{x \in \R^d} |f(x)|
\]
is its supremum norm, and we set
\[
\|f\|_{\infty,A}
= \sup_{x \in A} |f(x)|
\]
for $A \subseteq \Rd$.

A finite collection $f_1, \dots, f_N:\Rd \rightarrow \R$
  is called an $L_p$ $\varepsilon$--covering of $\F$ on $x_1^n$
  if for all $f \in \F$
  \[
  \min_{1 \leq j\leq N}
  \left(
  \frac{1}{n} \sum_{k=1}^n |f(x_k)-f_j(x_k)|^p
  \right)^{1/p} \leq \varepsilon
  \]
  hold.
  The $L_p$ $\varepsilon$--covering number of $\F$ on $x_1^n$
  is the  size $N$ of the smallest $L_p$ $\varepsilon$--covering
  of $\F$ on $x_1^n$ and is denoted by $\Nu_p(\varepsilon,\F,x_1^n)$.

For $z \in \R$ and $\beta>0$ we define
$T_\beta z = \max\{-\beta, \min\{\beta,z\}\}$. If $f:\R^d \rightarrow
\R$
is a function  then we set
$
(T_{\beta} f)(x)=
T_{\beta} \left( f(x) \right)$.

\subsection{Outline}
\label{se1sub6}
The main result is formulated in Section 2 and proven
in Section 3.

\section{Estimation of a $(p,C)$--smooth regression function}
\label{se2}

Throughout the paper we let
$\sigma(x)=1/(1+e^{-x})$
be the logistic squasher. We define the topology of our
neural networks as follows: We
let $K_n, L, r \in \N$ be parameters of our estimate and using
these parameters we
 set
\begin{equation}\label{se2eq1}
f_\bw(x) = \sum_{j=1}^{K_n} w_{1,1,j}^{(L)} \cdot f_{j,1}^{(L)}(x)
\end{equation}
for some $w_{1,1,1}^{(L)}, \dots, w_{1,1,K_n}^{(L)} \in \mathbb{R}$, where
$f_{j,1}^{(L)}=f_{\bw,j,1}^{(L)}$ are recursively defined by
\begin{equation}
  \label{se2eq2}
f_{k,i}^{(l)}(x) = 
f_{\bw,k,i}^{(l)}(x) = 
\sigma\left(\sum_{j=1}^{r} w_{k,i,j}^{(l-1)}\cdot f_{k,j}^{(l-1)}(x) + w_{k,i,0}^{(l-1)} \right)
\end{equation}
for some $w_{k,i,0}^{(l-1)}, \dots, w_{k,i, r}^{(l-1)} \in \mathbb{R}$
$(l=2, \dots, L)$
and
\begin{equation}
  \label{se2eq3}
f_{k,i}^{(1)}(x) = 
f_{\bw,k,i}^{(1)}(x) = 
\sigma \left(\sum_{j=1}^d w_{k,i,j}^{(0)}\cdot x^{(j)} + w_{k,i,0}^{(0)} \right)
\end{equation}
for some $w_{k,i,0}^{(0)}, \dots, w_{k,i,d}^{(0)} \in \mathbb{R}$.

This means that we consider neural networks which consist of $K_n$ fully
connected
neural networks of depth $L$ and width $r$ computed in parallel and compute
 a linear combination of the outputs of these $K_n$ neural
networks.
The weights in the $k$-th such network are denoted by
$(w_{k,i,j}^{(l)})_{i,j,l}$, where
$w_{k,i,j}^{(l)}$ is the weight between neuron $j$ in layer
$l$ and neuron $i$ in layer $l+1$.

We initialize the weights $\bw^{(0)}=((\bw^{(0)})_{k,i,j}^{(l))})_{k,i,j,l}$ as
follows: We set
\begin{equation}
\label{se2eq4}
(\bw^{(0)})_{1,1,k}^{(L)}=0
\quad (k=1, \dots, K_n),
\end{equation}
 we choose $(\bw^{(0)})_{k,i,j}^{(l)}$ uniformly distributed on
$[-c_1, c_1]$ if $l \in \{1, \dots, L-1\}$, and we
choose
$(\bw^{(0)})_{k,i,j}^{(0)}$ uniformly distributed on
$[-c_2 \cdot (\log n) \cdot n^\tau, c_2 \cdot (\log n) \cdot n^\tau]$, where $c_1, c_2, \tau>0$ are parameters of the estimate.
Here the random values are defined such that all components
of $\bw^{(0)}$ are independent.

After initialization of the weights we perform $t_n \in \N$ gradient
descent
steps each with a step size $\lambda_n>0$. Here we try to minimize
the empirical $L_2$ risk
\begin{equation}
\label{se2eq5}
F_n(\bw)
= \frac{1}{n} \sum_{i=1}^n | Y_i - f_\bw (X_i)|^2.
\end{equation}
 To do this we set
\begin{equation}
\label{se2eq6}
\bw^{(t)}
=
\bw^{(t-1)}
-
\lambda_n \cdot \nabla_{\bw} F_n(\bw^{(t-1)})
\quad
(t=1, \dots, t_n).
\end{equation}
Finally we define our estimate as a truncated version of the neural
network with weight vector $\bw^{(t_n)}$, i.e., we set
\begin{equation}
\label{se2eq7}
m_n(x)= T_{\beta_n} (f_{\bw^{(t_n)}}(x))
\end{equation}
 where $\beta_n = c_3 \cdot \log n$ and $T_{\beta} z
= \max\{ \min\{z, \beta\}, - \beta\}$ for $z \in \R$
and $\beta>0$.

Our main result is the following bound on the expected $L_2$
error of this estimate.

\begin{theorem}
  \label{th1}
Let $n \in \N$,
let $(X,Y)$, $(X_1,Y_n)$, \dots, $(X_n,Y_n)$
be independent and identically distributed $\Rd \times \R$--valued random variables such that $supp(X)$ is bounded and that
\begin{equation}
\label{th1eq1}
\EXP\left\{
e^{c_4 \cdot Y^2}
\right\}
< \infty
\end{equation}
holds for some $c_4>0$.
Let $p,C>0$ where $p=q+\beta$ for some $q \in \N_0$
and $\beta \in (0,1]$ with $p \geq 1/2$,
  and assume that the regression function
$m:\R^d \rightarrow \R$ is $(p,C)$--smooth.

  Set $\beta_n=c_3 \cdot \log n$ for some $c_3>0$ which satisfies
  $c_3 \cdot c_4 \geq 2$. Let
  $K_n \in \N$ be such that for some $\kappa>0$
  \[
  \frac{K_n}{n^\kappa} \rightarrow 0 \quad (n \rightarrow \infty)
  \quad \mbox{and} \quad
  \frac{K_n}{n^{
4 \cdot r \cdot (r+1) \cdot (L-1)  + r \cdot (4d+6) +6
  }}
  \rightarrow \infty \quad (n \rightarrow \infty).
  \]
  Set
\[
L=\lceil \log_2(q+d) \rceil+1, \quad r=2 \cdot \lceil (2p+d)^2 \rceil,
\quad
\tau=\frac{1}{2p+d}, \quad
\lambda_n=\frac{c_5}{ n \cdot K_n^3}
\]
and
\[
t_n=\left\lceil
c_6 \cdot \frac{K_n^3}{\beta_n}
\right\rceil
\]
for some $c_5, c_6 >0$.
Let $\sigma(x)=1/(1+e^{-x})$ be the logistic squasher,
let $c_1, c_2 >0$ be sufficiently large,
and define the estimate $m_n$ as above.

Then we have for any $\epsilon>0$:
\[
\EXP \int | m_n(x)-m(x)|^2 \PROB_X (dx)
\leq c_7 \cdot n^{- \frac{2p}{2p+d} + \epsilon}.
\]
  \end{theorem}

\noindent
    {\bf Remark 1.} By combining the approximation result
    derived in the proof of Theorem \ref{th1} with the proof
    strategy presented in Kohler and Krzy\.zak (2022) and Drews
    and Kohler (2023) it is possible to show that the rate
    of convergence in Theorem \ref{th1} also  holds if
    the inner weights are not learned
at all and gradient descent is applied only to the weights in the
output level.

\noindent
    {\bf Remark 2.} It should be easy to extend the above result
    to interaction models as in Kohler and Krzy\.zak (2022) and Drews
    and Kohler (2023), i.e., to modify the estimate in Theorem \ref{th1}
    such that it achieves the rate of convergence
    \[
    n^{- \frac{2p}{2p+d^*} + \epsilon}
    \]
    in case that the regression function is given by a $(p,C)$--smooth
    interaction model where each function in the sum depends on at
    most $d^* \in \{1, \dots, d\}$ of the $d$ components of $X$.

    \noindent
        {\bf Remark 3.}
        It is an open problem whether the above result can be extended
        to the case of an hierarchical composition model.

    \section{Proof of Theorem \ref{th1}}
\label{se4}
Before we present the proof of Theorem \ref{th1} we present
in separate subsections the key auxiliary results needed in the
proof concerning optimization, approximation and generalization.

\subsection{Neural network optimization}

Our first lemma is our main tool to analyze gradient descent.
In it we relate the gradient descent of our deep neural network to the
gradient descent of the linear Taylor polynomial of the deep network,
and use methods for the analysis of gradient descent applied to
smooth convex functions in order to analyze the latter.

\begin{lemma}
  \label{le1}
  Let $d, J_n \in \N$, and for $\bw \in \R^{J_n}$ let
$
  f_{\bw}: \R^d \rightarrow \R
  $
  be a (deep) neural network with weight vector $\bw$.
  Assume that for each $x \in \R^d$
  \[
\bw \mapsto f_\bw(x)
  \]
is a continuously differentiable function on $\R^{J_n}$.
  Let
\[
F_n(\bw) = \frac{1}{n} \sum_{i=1}^n |Y_i - f_{\bw}(X_i)|^2
\]
be the empirical $L_2$ risk of $f_{\bw}$, and use gradient
descent in order to minimize $F_n(\bw)$. To do this, choose
a starting weight vector $\bw^{(0)} \in \R^{J_n}$, choose
$\delta_n \geq 0$ and let
\[
A \subset \left\{
\bw \in \R^{J_n} \, : \, \|\bw-\bw^{(0)}\| \leq \delta_n
\right\}
\]
be a closed and convex set of weight vectors. Choose a stepsize
$\lambda_n \geq 0$ and a number of gradient descent steps $t_n \in \N$
and compute
\[
\bw^{(t+1)} = Proj_A \left(
\bw^{(t)}
-
\lambda_n
\cdot \nabla_{\bw} F_n( \bw^{(t)} )
\right)
\]
for $t=0, \dots, t_n-1$.

Let $C_n, D_n \geq 0$, $\beta_n \geq 1$ and assume
\begin{equation}
  \label{le1eq1}
  \sum_{j=1}^{J_n}
  \left|
  \frac{\partial}{\partial w^{(j)}}
  f_{\bw_1}(x)
-
\frac{\partial}{\partial w^{(j)}}
f_{\bw_2}(x)
\right|^2
   \leq C_n^2 \cdot
  \| \bw_1 - \bw_2 \|^2
  \end{equation}
  for all  $\bw_1, \bw_2 \in A$, $x \in \{X_1, \dots, X_n\}$,
\begin{equation}
  \label{le1eq2}
  \| \nabla_\bw F_n(\bw) \| \leq D_n
  \quad
\mbox{for all }  \bw \in A,  
  \end{equation}
\begin{equation}
  \label{le1eq3}
  |Y_i| \leq \beta_n \quad (i=1, \dots, n) 
\end{equation}
and
\begin{equation}
  \label{le1eq4}
C_n \cdot \delta_n^2 \leq 1.
  \end{equation}
Let $\bw^* \in A$ and assume
\begin{equation}
  \label{le1eq5}
  |f_{\bw^*}(x)| \leq \beta_n \quad (x \in \{X_1, \dots, X_n\}).
\end{equation}
Then  
\begin{eqnarray*}
  &&
  \min_{t=0, \dots, t_n-1} F_n(\bw^{(t)})
  \leq
  F_n(\bw^*) + \frac{\| \bw^* - \bw^{(0)} \|^2}{2 \cdot \lambda_n \cdot t_n}
  +
  12 \cdot \beta_n \cdot C_n \cdot \delta_n^2 +
  \frac{1}{2} \cdot \lambda_n \cdot D_n^2.
  \end{eqnarray*}
  \end{lemma}

\noindent
    {\bf Proof.}
    The basic idea of the proof is to analyze the
    gradient descent by relating it to the gradient
    descent of the linear Taylor polynomial of $f_\bw$. To do this,
    we define
    for $\bw_0, \bw \in \R^{J_n}$ the linear Taylor polynomial
    of $f_\bw(x)$ around $\bw_0$ by
    \[
    f_{lin, \bw_0, \bw}(x)
    =
    f_{\bw_0}(x)+
    \sum_{j=1}^{J_n}
    \frac{\partial f_{\bw_0}(x)}{ \partial \bw^{(j)}}
      \cdot
      (\bw^{(j)} - \bw_0^{(j)})
    \]
    and introduce the empirical $L_2$ risk of this
    linear approximation of $f_\bw$ by
    \[
    F_{n,lin,\bw_0}(\bw)
    =
    \frac{1}{n} \sum_{i=1}^n |Y_i - f_{lin,\bw_0,\bw}(X_i)|^2.
    \]
Let $\alpha \in [0,1]$ and $\bw_1, \bw_2 \in \R^{J_n}$. Then
\begin{eqnarray*}
  &&
  f_{lin, \bw_0, \alpha \cdot \bw_1 + (1-\alpha) \cdot \bw_2 }(x)
  \\
  &&
    =
    f_{\bw_0}(x)+
    \sum_{j=1}^{J_n}
    \frac{\partial f_{\bw_0}(x)}{ \partial \bw^{(j)}}
      \cdot
      (\alpha \cdot \bw_1^{(j)} + (1-\alpha) \cdot \bw_2^{(j)} - \bw_0^{(j)})
      \\
      &&
      =\alpha \cdot f_{\bw_0}(x)+(1-\alpha) \cdot f_{\bw_0}(x)+
      \alpha \cdot
    \sum_{j=1}^{J_n}
    \frac{\partial f_{\bw_0}(x)}{ \partial \bw^{(j)}}
      \cdot
      (\bw_1^{(j)} - \bw_0^{(j)})
      \\
      &&
      \quad
      +
      (1-\alpha) \cdot
    \sum_{j=1}^{J_n}
    \frac{\partial f_{\bw_0}(x)}{ \partial \bw^{(j)}}
      \cdot
      (\bw_1^{(j)} - \bw_0^{(j)})
      \\
      &&
      =
      \alpha \cdot f_{lin, \bw_0, \bw_1 }(x)
      +
      (1-\alpha) \cdot f_{lin, \bw_0, \bw_2 }(x),
  \end{eqnarray*}
which implies
\begin{eqnarray*}
  &&
  F_{n,lin,\bw_0}(\alpha \cdot \bw_1 + (1-\alpha) \cdot \bw_2 )
  \\
  &&
    =
    \frac{1}{n} \sum_{i=1}^n |\alpha \cdot (Y_i - f_{lin,\bw_0,\bw_1}(X_i))
    + (1-\alpha) \cdot (Y_i - f_{lin,\bw_0,\bw_2}(X_i))|^2
    \\
    &&
    \leq
    \alpha \cdot   F_{n,lin,\bw_0}(\bw_1) + (1-\alpha) \cdot   F_{n,lin,\bw_0}(\bw_2).
  \end{eqnarray*}
Hence $F_{n,lin,\bw_0}(\bw)$ is as a function of $\bw$ a convex function.

Because of 
$   f_{lin, \bw_0, \bw_0}(x)
    =
    f_{\bw_0}(x)$ and
    $   \nabla_\bw f_{lin, \bw_0, \bw_0}(x)
    =
    \nabla_w f_{\bw_0}(x)$ we have
    \[
    F_{n,lin,\bw^{(t)}}(\bw^{(t)}) = F_n(\bw^{(t)})
      \quad \mbox{and} \quad
    \nabla_w F_{n,lin,\bw^{(t)}}(\bw^{(t)})
    =
    \nabla_w F_{n}(\bw^{(t)}),
    \]
    hence $\bw^{(t+1)}$ is computed from $\bw^{(t)}$ by one gradient descent
    step
    \[
\bw^{(t+1)} = Proj_A \left(
\bw^{(t)}
-
\lambda_n
\cdot \nabla_{\bw} F_{n,lin,\bw^{(t)}}( \bw^{(t)} )
\right)
    \]
    applied to the convex function $ F_{n,lin,\bw^{(t)}}( \bw )$.
    This will enable us to use techniques for the analysis of the gradient descent
    for convex functions in order to analyze the gradient descent applied
    to the nonconvex function $F_n(\bw)$.

    In order to do this we observe
    \begin{eqnarray*}
      &&
  \min_{t=0, \dots, t_n-1} F_n(\bw^{(t)})
  -
  F_n(\bw^*)
  \\
  &&
  \leq
  \frac{1}{t_n} \sum_{t=0}^{t_n-1}
  (  F_n(\bw^{(t)}) -
  F_n(\bw^*))
  \\
  &&
  =
  \frac{1}{t_n} \sum_{t=0}^{t_n-1}
  (  F_{n,lin, \bw^{(t)}}(\bw^{(t)}) -
  F_{n,lin,\bw^{(t)}}(\bw^{*}))
  +
  \frac{1}{t_n} \sum_{t=0}^{t_n-1}
      ( F_{n,lin,\bw^{(t)}}(\bw^{*})-F_n(\bw^*))
      \\
      &&
      =: T_{1,n} + T_{2,n}.
    \end{eqnarray*}
    Next we show that
    assumption (\ref{le1eq1}) implies
    \[
|    f_{\bw}(x)-f_{lin,\bw_0,\bw}(x)| \leq  \frac{1}{2} \cdot C_n \cdot \| \bw - \bw_0\|^2
\]
for all $x \in \{X_1, \dots, X_n\}$ and all $\bw_0, \bw \in A$.
To do this, set
\[
H(s)=f_{\bw_0 + s \cdot (\bw-\bw_0)}(x) \quad \mbox{for } s \in [0,1].
\]
Let 
$\bw, \bw_0 \in A$.
Then $A$ convex implies
\[
\bw_0+s \cdot (\bw-\bw_0)
=
(1-s) \cdot \bw_0+s \cdot \bw \in A
\]
for all $s \in [0,1]$, hence we can conclude from (\ref{le1eq1})
\begin{eqnarray*}
  &&
  |    f_{\bw}(x)-f_{lin,\bw_0,\bw}(x)|
  \\
  &&
  =
  |  f_{\bw}(x)-
    f_{\bw_0}(x)-
    \sum_{j=1}^{J_n}
    \frac{\partial f_{\bw_0}(x)}{ \partial \bw^{(j)}}
      \cdot
      (\bw^{(j)} - \bw_0^{(j)})
      |
      \\
      &&
      =
      |H(1)-H(0) - \sum_{j=1}^{J_n}
    \frac{\partial f_{\bw_0}(x)}{ \partial \bw^{(j)}}
      \cdot
      (\bw^{(j)} - \bw_0^{(j)})
      |
      \\
      &&
      =
      |
      \int_0^1 H^\prime(s) \, ds
      - \sum_{j=1}^{J_n}
    \frac{\partial f_{\bw_0}(x)}{ \partial \bw^{(j)}}
      \cdot
      (\bw^{(j)} - \bw_0^{(j)})
      |
      \\
      &&
      =
      |
      \int_0^1
      \sum_{j=1}^{J_n}
    \frac{\partial f_{\bw_0+s \cdot (\bw-\bw_0)}(x)}{ \partial \bw^{(j)}}
      \cdot
      (\bw^{(j)} - \bw_0^{(j)})
      \, ds
      - \sum_{j=1}^{J_n}
    \frac{\partial f_{\bw_0}(x)}{ \partial \bw^{(j)}}
      \cdot
      (\bw^{(j)} - \bw_0^{(j)})
      |
      \\
      &&
      =
      |
      \int_0^1
      \sum_{j=1}^{J_n}
      (    \frac{\partial f_{\bw_0+s \cdot (\bw-\bw_0)}(x)}{ \partial \bw^{(j)}}
      -
    \frac{\partial f_{\bw_0}(x)}{ \partial \bw^{(j)}})  
      \cdot
      (\bw^{(j)} - \bw_0^{(j)})
      \, ds
      |
      \\
      &&
      \leq
      \int_0^1
      \sum_{j=1}^{J_n}
      |    \frac{\partial f_{\bw_0+s \cdot (\bw-\bw_0)}(x)}{ \partial \bw^{(j)}}
      -
    \frac{\partial f_{\bw_0}(x)}{ \partial \bw^{(j)}}|  
      \cdot
      |\bw^{(j)} - \bw_0^{(j)}|
      \, ds
      \\
      &&
      \leq
        \int_0^1
        \sqrt{
 \sum_{j=1}^{J_n}
      |    \frac{\partial f_{\bw_0+s \cdot (\bw-\bw_0)}(x)}{ \partial \bw^{(j)}}
      -
    \frac{\partial f_{\bw_0}(x)}{ \partial \bw^{(j)}}|^2 
        }
        \cdot
        \| \bw - \bw_0\|  \, ds
        \\
        &&
        \leq
        \int_0^1
        \sqrt{
          C_n^2 \cdot \| \bw_0+s \cdot (\bw-\bw_0) - \bw_0\|^2}
         \cdot
         \| \bw - \bw_0\|  \, ds
         \\
         &&
         \leq
         C_n \cdot  \| \bw - \bw_0\|^2 \cdot \int_0^1 s \, ds
         =
         \frac{1}{2} \cdot  C_n \cdot  \| \bw - \bw_0\|^2.
  \end{eqnarray*}

Using (\ref{le1eq3})--(\ref{le1eq5}) we can conclude
for all $\bw_0\in A$
\begin{eqnarray*}
  &&
|  F_n(\bw^*)
  -
    F_{n,lin,\bw_0}(\bw^*)
    |
    \\
    &&
    \leq
    \frac{1}{n} \sum_{i=1}^n |Y_i - f_{\bw^*}(X_i) + Y_i - f_{lin,\bw_0,\bw^*}(X_i)|
    \cdot
    |f_{\bw^*}(X_i)  - f_{lin,\bw_0,\bw^*}(X_i)|
    \\
    &&
    \leq
    \frac{1}{n} \sum_{i=1}^n (4 \cdot \beta_n +   
    \frac{1}{2} \cdot C_n \cdot \| \bw^* - \bw_0\|^2) \cdot
    \frac{1}{2} \cdot C_n \cdot \| \bw^* - \bw_0\|^2
    \\
    &&
    \leq
    \frac{1}{n} \sum_{i=1}^n (4 \cdot \beta_n + \frac{1}{2} \cdot  C_n \cdot
    4 \delta_n^2) \cdot \frac{1}{2} \cdot C_n \cdot 4 \delta_n^2
    \\
    &&
    \leq
    12 \cdot \beta_n \cdot C_n \cdot \delta_n^2.
  \end{eqnarray*}
This proves
\[
T_{2,n}
   \leq
    12 \cdot \beta_n \cdot C_n \cdot \delta_n^2,
\]
and it suffices to show
\begin{equation}
  \label{le1peq1}
T_{1,n}
\leq
\frac{\| \bw^* - \bw^{(0)} \|^2}{2 \cdot \lambda_n \cdot t_n}
 + \frac{1}{2} \cdot \lambda_n \cdot D_n^2.
\end{equation}

    The convexity of $ F_{n,lin,\bw^{(t)}}( \bw )$ together with $\bw^* \in A$
    implies
    \begin{eqnarray*}
      &&
      F_{n,lin,\bw^{(t)}}(\bw^{(t)}) - F_{n,lin,\bw^{(t)}}(\bw^{*})
      \\
      &&
      \leq
\;      < \nabla_\bw       F_{n,lin,\bw^{(t)}}(\bw^{(t)}), \bw^{(t)}-\bw^* >
      \\
      &&
      =
  \;          < \nabla_\bw       F_{n}(\bw^{(t)}), \bw^{(t)}-\bw^* >
            \\
            &&
            =
            \frac{1}{2 \cdot \lambda_n}
            \cdot 2 \cdot < \lambda_n  \cdot \nabla_\bw    F_{n}(\bw^{(t)}), \bw^{(t)}-\bw^* >
            \\
            &&
            =
             \frac{1}{2 \cdot \lambda_n}
             \cdot
             \left(
\| \bw^{(t)}-\bw^*\|^2 - \|\bw^{(t)}-\bw^* - \lambda_n  \cdot \nabla_\bw    F_{n}(\bw^{(t)})\|^2 + \| \lambda_n  \cdot \nabla_\bw    F_{n}(\bw^{(t)})\|^2
\right)
\\
&&
     =
             \frac{1}{2 \cdot \lambda_n}
             \cdot
             \left(
\| \bw^{(t)}-\bw^*\|^2 - \|\bw^{(t)} - \lambda_n  \cdot \nabla_\bw    F_{n}(\bw^{(t)})- \bw^* \|^2 
\right)
+
\frac{1}{2} \cdot \lambda_n \cdot
\| \nabla_\bw    F_{n}(\bw^{(t)})\|^2
\\
&&
\leq
             \frac{1}{2 \cdot \lambda_n}
             \cdot
             \left(
             \| \bw^{(t)}-\bw^*\|^2 - \|Proj_A \left(
             \bw^{(t)} - \lambda_n  \cdot \nabla_\bw
             F_{n}(\bw^{(t)}) \right)- \bw^* \|^2 
\right)
\\
&&
\quad
+
\frac{1}{2} \cdot \lambda_n \cdot
\| \nabla_\bw    F_{n}(\bw^{(t)})\|^2
\\
&&
=      \frac{1}{2 \cdot \lambda_n}
             \cdot
             \left(
             \| \bw^{(t)}-\bw^*\|^2 -
             \| \bw^{(t+1)}-\bw^*\|^2 \right)
             +
\frac{1}{2} \cdot \lambda_n \cdot
\| \nabla_\bw    F_{n}(\bw^{(t)})\|^2.
      \end{eqnarray*}
    This together with (\ref{le1eq2}) implies
    \begin{eqnarray*}
      T_{1,n}
      & \leq &
      \frac{1}{t_n} \sum_{t=0}^{t_n-1}
      \left(
 \frac{1}{2 \cdot \lambda_n}
             \cdot
             \left(
             \| \bw^{(t)}-\bw^*\|^2 -
             \| \bw^{(t+1)}-\bw^*\|^2 \right)
             +
\frac{1}{2} \cdot \lambda_n \cdot
\| \nabla_\bw    F_{n}(\bw^{(t)})\|^2
\right)
\\
&\leq&
\frac{
\| \bw^{(0)}-\bw^*\|^2 
}{2 \cdot \lambda_n \cdot t_n}
+ \frac{1}{2} \cdot \lambda_n \cdot D_n^2,
      \end{eqnarray*}
which proves (\ref{le1peq1}). \hfill $\Box$%

Next we consider the topology of the deep neural
network introduced in Section \ref{se2}
(cf., (\ref{se2eq1})-(\ref{se2eq3}))
and investigate when the assumptions of Lemma \ref{le1}
are satisfied.
  
Our next lemma considers inequality (\ref{le1eq1}) in this case.

\begin{lemma}
  \label{le2}
  Let $\sigma$ be the logistic squasher.
  Let $a,B_n, \gamma_n^* \geq 1$, $L,r \in \N$ and define the deep neural network
  $f_\bw:\R^d \rightarrow \R$ with weight vector $\bw$
  by (\ref{se2eq1})--(\ref{se2eq3}). Assume that the weight
  vectors $\bw_1$ and $\bw_2$ satisfy
  \[
  |w_{1,1,k}^{(L)}| \leq \gamma_n^*
  \quad \mbox{and} \quad
  |w_{k,i,j}^{(l)}| \leq B_n 
  \]
  for all $l \in \{1, \dots, L-1\}$.
  Then we have for any $x \in [-a,a]^d$
  \[
  \sum_{k,i,j,l}
  \left|
  \frac{\partial}{\partial w_{k,i,j}^{(l)}}
  f_{\bw_1}(x)
-
\frac{\partial}{\partial w_{k,i,j}^{(l)}}
f_{\bw_2}(x)
\right|^2
   \leq c_8  \cdot B_n^{4L} \cdot (\gamma_n^*)^2  \cdot
  \| \bw_1 - \bw_2 \|^2
  \]
  for some $c_8=c_8(d,L,r,a)>0$.
\end{lemma}

\noindent
    {\bf Proof.}
    We have
    \begin{eqnarray*}
      &&
  \sum_{k,i,j,l}
  \left|
  \frac{\partial}{\partial w_{k,i,j}^{(l)}}
  f_{\bw_1}(x)
-
\frac{\partial}{\partial w_{k,i,j}^{(l)}}
f_{\bw_2}(x)
\right|^2
\\
&&
=
\sum_{k=1}^{K_n}
|f_{\bw_1,k,1}^{(L)}(x) - f_{\bw_2, k,1}^{(L)}(x)|^2
\\
&&
\quad
+
\sum_{k=1}^{K_n}
\sum_{i,j,l : l<L}
\left|
(\bw_1)_{1,1,k}^{(L)} \cdot \frac{\partial}{\partial \bw_{k,i,j}^{(l)}}
  f_{\bw_1,k,1}^{(L)}(x)
  -
(\bw_2)_{1,1,k}^{(L)} \cdot \frac{\partial}{\partial \bw_{k,i,j}^{(l)}}
  f_{\bw_2,k,1}^{(L)}(x)  
  \right|^2\\
  &&
  \leq
\sum_{k=1}^{K_n}
|f_{\bw_1,k,1}^{(L)}(x) - f_{\bw_2, k,1}^{(L)}(x)|^2
\\
&&
\quad
+ 2 \cdot
\sum_{k=1}^{K_n}
\sum_{i,j,l: l<L}
\left|
(\bw_1)_{1,1,k}^{(L)}
  -
(\bw_2)_{1,1,k}^{(L)}  
  \right|^2
\cdot |  \frac{\partial}{\partial \bw_{k,i,j}^{(l)}} f_{\bw_1,k,1}^{(L)}(x)|^2
  \\
    &&
\quad
+ 2 \cdot
\sum_{k=1}^{K_n}
\sum_{i,j,l: l<L}
\left| (\bw_2)_{1,1,k}^{(L)} \right|^2 \cdot
\left|
\frac{\partial}{\partial \bw_{k,i,j}^{(l)}}
  f_{\bw_1,k,1}^{(L)}(x)
  -
 \frac{\partial}{\partial \bw_{k,i,j}^{(l)}}
  f_{\bw_2,k,1}^{(L)}(x)  
  \right|^2.
    \end{eqnarray*}
     The chain rule implies
    \begin{eqnarray}
      &&
      \frac{\partial f_{\bw,k,1}^{(L)}}{\partial w_{k,i,j}^{(l)}}(x)
        =
  \sum_{s_{l+2}=1}^{r} \dots \sum_{s_{L-1}=1}^{r}
  f_{k,j}^{(l)}(x)
  \cdot
\sigma^\prime \left(\sum_{t=1}^{r} w_{k,i,t}^{(l)} \cdot f_{k,t}^{(l)}(x) + w_{k,i,0}^{(l)} \right)
  \nonumber \\
  && \quad
  \cdot
  w_{k,s_{l+2},i}^{(l+1)} \cdot
\sigma^\prime \left(\sum_{t=1}^{r} w_{k,s_{l+2},t}^{(l+1)} \cdot f_{k,t}^{(l+1)}(x) + w_{k,s_{l+2},0}^{(l+1)} \right)
  \cdot
  w_{k,s_{l+3},s_{l+2}}^{(l+2)}
   \nonumber \\
  && \quad
  \cdot
\sigma^\prime \left(\sum_{t=1}^{r} w_{k,s_{l+3},t}^{(l+2)} \cdot f_{k,t}^{(l+2)}(x) + w_{k,s_{l+3},0}^{(l+2)} \right)
  \cdots
  w_{k,s_{L-1},s_{L-2}}^{(L-2)}
 \nonumber \\
  && \quad
  \cdot
\sigma^\prime \left(\sum_{t=1}^{r} w_{k,s_{L-1},t}^{(L-2)} \cdot f_{k,t}^{(L-2)}(x) + w_{k,s_{L-1},0}^{(L-2)} \right)
   \cdot
   w_{k,1,s_{L-1}}^{(L-1)}
 \nonumber \\
  && \quad
   \cdot
\sigma^\prime \left(\sum_{t=1}^{r} w_{k,1,t}^{(L-1)} \cdot f_{k,t}^{(L-1)}(x) + w_{k,1,0}^{(L-1)} \right),
     \label{ple2eq1}
      \end{eqnarray}
where we have used the abbreviations
\[
f_{k,j}^{(0)}(x)
=
\left\{
\begin{array}{ll}
  x^{(j)} & \mbox{if } j \in \{1,\dots,d\} \\
  1 & \mbox{if } j=0
\end{array}
\right.
\]
and
\[
f_{k,0}^{(l)}(x)=1 \quad (l=1, \dots, L-1).
\]
If $f_{i,1}$, \dots, $f_{i,L}$ are real--valued functions defined on $\R^{J_n}$
where $f_{i,l}$ is bounded in absolute value by $B_{i,l} \geq 1$ and Lipschitz
continuous (w.r.t. $\|\cdot\|_\infty$) with Lipschitz constant $C_{i,l} \geq 1$ $(i=1, \dots, r)$
then
\begin{eqnarray*}
  &&
| \sum_{i=1}^r \prod_{l=1}^L f_{i,l}(\bw_1)
-
\sum_{i=1}^r \prod_{l=1}^L f_{i,l}(\bw_2)|
\\
&&
\leq
\sum_{i=1}^r
\sum_{j=1}^L
\prod_{l=1}^{j-1} |f_{i,l}(\bw_1)|
\cdot
|f_{j,l}(\bw_1) - f_{j,l}(\bw_2)|
\prod_{l=j+1}^L |f_{i,l}(\bw_2)|
\\
&&
\leq
r \cdot L \cdot \max_{i=1, \dots, r}  \prod_{l=1}^L B_{i,l}
\cdot
\max_{i=1, \dots, r} \max_{l=1, \dots, L} C_{i,l}
\cdot \|\bw_1-\bw_2\|_\infty.
\end{eqnarray*}
Using this,
\[
0 \leq \sigma(x) \leq 1 \quad \mbox{and} \quad
|\sigma^\prime(x)|=|\sigma(x) \cdot (1-\sigma(x))| \leq 1
\]
and 
\begin{eqnarray*}
&&|  f_{\bw_1,k,j}^{(l)}(x)-  f_{\bw_2,k,j}^{(l)}(x)|\\
&& \leq c_9 \cdot a \cdot (\max\{ 2r, d\}+1)^l \cdot B_n^{l-1}  \cdot
\| ((\bw_1)_{k,\bar{i},\bar{j}}^{(\bar{l})})_{\bar{i},\bar{j},\bar{l}} -
((\bw_2)_{k,\bar{i},\bar{j}}^{(\bar{l})})_{\bar{i},\bar{j},\bar{l}}\|_\infty\\
&& \leq c_9 \cdot a \cdot (\max\{2r,d\}+1)^l \cdot B_n^{l-1}  \cdot
\| ((\bw_1)_{k,\bar{i},\bar{j}}^{(\bar{l})})_{\bar{i},\bar{j},\bar{l}} -
((\bw_2)_{k,\bar{i},\bar{j}}^{(\bar{l})})_{\bar{i},\bar{j},\bar{l}}\|
\end{eqnarray*}
(which can be easily shown by induction on $l$)
we get 
    \begin{eqnarray*}
      &&
  \sum_{k,i,j,l}
  \left|
  \frac{\partial}{\partial w_{k,i,j}^{(l)}}
  f_{\bw_1}(x)
-
\frac{\partial}{\partial w_{k,i,j}^{(l)}}
f_{\bw_2}(x)
\right|^2
\\
&&
\leq
c_{10} \cdot a^2 \cdot (2r+d)^{2L} \cdot B_n^{2L}  \cdot \|\bw_1-\bw_2\|^2
\\
&&
\quad
+c_{11} \cdot L \cdot
(r \cdot (r+d)) \cdot r^{2L} \cdot a^2 \cdot B_n^{2L}  \cdot \|\bw_1-\bw_2\|^2
\\
&&
\quad
+c_{12} \cdot (\gamma_n^*)^2 \cdot L \cdot (r \cdot (r+d))
\cdot r^{2L} \cdot (3L)^2 \cdot a^4 \cdot B_n^{4L} 
\cdot (2r+d)^{2L} \cdot \|\bw_1-\bw_2\|^2
\\
&&
\leq c_{13} \cdot (\gamma_n^*)^2 \cdot L^3 \cdot (2r+d)^{4L+2}
\cdot B_n^{4L} \cdot a^4 \cdot \|\bw_1-\bw_2\|^2.
\end{eqnarray*}
    \hfill $\Box$

    Next we consider inequality (\ref{le1eq2})
    in case of the special topology of our networks.

    \begin{lemma}
      \label{le3}
        Let $\sigma$ be the logistic squasher.
        Let $a,\beta_n, B_n, \gamma_n^* \geq 1$, $K_n,L,r \in \N$ and
      define the deep neural network
  $f_\bw:\R^d \rightarrow \R$ with weight vector $\bw$
      by (\ref{se2eq1})--(\ref{se2eq3}), and assume
      \[
      X_i \in [-a,a]^d
            \quad
  \mbox{and} \quad
      |Y_i| \leq \beta_n \quad (i=1, \dots, n)
\]
and
  \[
  |w_{1,1,k}^{(L)}| \leq \gamma_n^* \quad
  \mbox{and} \quad
  |w_{k,i,j}^{(l)}| \leq B_n 
  \]
  for all $l \in \{1, \dots, L-1\}$. Assume
  \[
K_n \cdot \gamma_n^* \geq \beta_n.
  \]

Then
\[
\| \nabla_\bw F_n(\bw) \| \leq c_{14} \cdot
K_n^{3/2} \cdot (\gamma_n^*)^2 \cdot B_n^L 
  \]
    for some $c_{14}=c_{14}(d,L,r,a)>0$.
      \end{lemma}

    \noindent
        {\bf Proof.}
        We have
        \[
|f_\bw(x)| \leq K_n \cdot \gamma_n^* ,
\]
which implies
        \begin{eqnarray*}
          &&
       \| \nabla_\bw F_n(\bw) \|^2    \\
      &&=
  \sum_{k,i,j,l}
  \left|
  \frac{1}{n}
  \sum_{s=1}^n 2 \cdot (Y_s - f_{\bw}(X_s)) \cdot  
  \frac{\partial}{\partial w_{k,i,j}^{(l)}}
  f_{\bw}(X_s) \cdot (-1)
  \right|^2
  \\
  && \leq
  8 \cdot (K_n \cdot \gamma_n^*)^2 \cdot
  \sum_{k,i,j,l}
  \max_{s=1, \dots, n}
  \left|
  \frac{\partial}{\partial w_{k,i,j}^{(l)}}
  f_{\bw}(X_s) 
  \right|^2
  \\
  && =
  8 \cdot (K_n \cdot \gamma_n^*)^2 \cdot \Bigg(
  \sum_{k=1}^{K_n}
  \max_{s=1, \dots, n}
|f_{\bw,k,1}^{(L)}(X_s)|^2
\\
&&
\quad \quad
+
\sum_{k=1}^{K_n}
\sum_{i,j,l : l<L}
  \max_{s=1, \dots, n}
\left|
\bw_{1,1,k}^{(L)} \cdot \frac{\partial}{\partial \bw_{k,i,j}^{(l)}}
  f_{\bw,k,1}^{(L)}(X_s)
  \right|^2
  \Bigg)
  \\
  &&
  \leq
  8 \cdot (K_n \cdot \gamma_n^*) \cdot \Bigg(
  K_n \cdot 1 \\
  && \quad + K_n \cdot L \cdot (r \cdot (r+d))^L \cdot (\gamma_n^*)^2
  \cdot
  \max_{s=1, \dots, n}
  \max_{k,i,j,l:l<L} \left|
\frac{\partial}{\partial \bw_{k,i,j}^{(l)}}
  f_{\bw,k,1}^{(L)}(X_s)
  \right|^2
  \Bigg)
  \\
  &&
  \leq
   8 \cdot (K_n \cdot \gamma_n^*)^2 \cdot \left(
  K_n \cdot 1 + K_n \cdot L \cdot (r \cdot (r+d)) \cdot (\gamma_n^*)^2
  \cdot  r^{2L}
  \cdot a^2 \cdot B_n^{2L} 
  \right),
\end{eqnarray*}        
        where the last inequality follows from (\ref{ple2eq1}), the assumptions
        on the weights and the bounds on the logistic squasher
        mentioned in the proof of Lemma \ref{le2}.
        \hfill $\Box$

        In order to avoid the projection step in Lemma \ref{le1}, we use
        the following localization lemma for gradient descent proven
        in Braun et al. (2023).

        \begin{lemma}
          \label{le4}
            Let
  $F:\R^K \rightarrow \R_+$
  be a nonnegative differentiable function.
  Let
  $t \in \N$, $L>0$, $\ba_0 \in \R^K$ and set
  \[
  \lambda=
\frac{1}{L}
\]
and
\[
\ba_{k+1}=\ba_k - \lambda \cdot (\nabla_{\ba} F)(\ba_k)
\quad
(k \in \{0,1, \dots, t-1\}).
\]
Assume
\begin{equation}
  \label{le4eq1}
  \left\|
 (\nabla_{\ba} F)(\ba)
  \right\|
  \leq
  \sqrt{
2 \cdot t \cdot L \cdot \max\{ F(\ba_0),1 \}
    }
\end{equation}
for all $\ba \in \R^K$ with
$\| \ba - \ba_0\| \leq \sqrt{2 \cdot t \cdot \max\{ F(\ba_0),1 \} / L}$,
and
\begin{equation}
  \label{le4eq2}
\left\|
(\nabla_{\ba} F)(\ba)
-
(\nabla_{\ba} F)(\bb)
  \right\|
  \leq
  L \cdot \|\ba - \bb \|
\end{equation}
for all $\ba, \bb \in \R^K$ satisfying
\begin{equation}
  \label{le4eq3}
  \| \ba - \ba_0\| \leq \sqrt{8 \cdot \frac{t}{L} \cdot \max\{ F(\ba_0),1 \}}
  \quad \mbox{and} \quad
  \| \bb - \ba_0\| \leq \sqrt{8 \cdot \frac{t}{L} \cdot \max\{ F(\ba_0),1 \}}.
\end{equation}
Then we have
\[
\|\ba_k-\ba_0\| \leq
\sqrt{
2 \cdot \frac{k}{L} \cdot (F(\ba_0)-F(\ba_k))
}
\quad
 \mbox{for all }
 k \in \{1, \dots,t\}
\]
 and
 \[
 F(\ba_k) \leq F(\ba_{k-1}) \quad
 \mbox{for all }
 k \in \{1, \dots,t\}.
 \]
          \end{lemma}

        \noindent
            {\bf Proof.} See Lemma A.1 in Braun et al. (2023) \hfill $\Box$

            Our next lemma helps us to verify the assumption (\ref{le4eq2}) of
            Lemma \ref{le4}.
            
    \begin{lemma}
      \label{le5}
      Let $\sigma: \R \rightarrow \R$ be the logistic squasher, let $f_\bw$
      be defined by (\ref{se2eq1})--(\ref{se2eq3}), and let
      $F_n$ be defined by (\ref{se2eq5}).
     Let $a \geq 1$,
     $\gamma_n^* \geq 1$, $B_n \geq 1$,  and assume
     $X_i \in [-a,a]^d$ $(i=1, \dots, n)$,
     \begin{equation}
     	\label{le5eq1}
     	\max\{ |(\bw_1)_{1,1,k}^{(L)}|, |(\bw_2)_{1,1,k}^{(L)}|\} \leq \gamma_n^* \quad (k=1,
     	\dots, K_n),
     \end{equation}
     \begin{equation}
     	\label{le5eq2}
     	\max\{|(\bw_1)_{k,i,j}^{(l)}|,|(\bw_2)_{k,i,j}^{(l)}|\} \leq B_n
     	\quad
     	\mbox{for } l=1, \dots, L-1
     \end{equation}
     and
  \[
K_n \cdot \gamma_n^* \geq \beta_n.
  \]
     Then we have
     \begin{eqnarray*}
     	&&
     	\| (\nabla_\bw F_n)(\bw_1) - (\nabla_\bw F_n)(\bw_2) \|
     	\leq
     	c_{15}   \cdot K_n^{3/2} \cdot B_n^{2L} \cdot (\gamma_n^*)^2 \cdot \|\bw_1-\bw_2\|.
     \end{eqnarray*}
\end{lemma}

\noindent
    {\bf Proof.}    
We have
    \begin{eqnarray*}
      &&
      \| \nabla_\bw F_n (\bw_1) -  \nabla_\bw F_n (\bw_2) \|^2
      \\
      &&
      =
      \sum_{k,i,j,l}
      \Bigg(
      \frac{2}{n}
      \sum_{s=1}^n
      (f_{\bw_1} (X_s) - Y_s) 
      \cdot
      \frac{\partial f_{\bw_1}}{\partial w_{k,i,j}^{(l)}}(X_s)
      -
      \sum_{k,i,j,l}
      \Bigg(
      \frac{2}{n}
      \sum_{s=1}^n
      (f_{\bw_2} (X_s)- Y_s) 
      \cdot
      \frac{\partial f_{\bw_2}}{\partial w_{k,i,j}^{(l)}}(X_s)
      \Bigg)^2
      \\
      &&
      \leq
      8 \cdot
      \sum_{k,i,j,l}
      \max_{s=1,\dots,n}
      \left(
      \frac{\partial f_{\bw_1}}{\partial w_{k,i,j}^{(l)}}(X_s)
      \right)^2 
    \cdot
      \frac{1}{n}
      \sum_{s=1}^n
      (f_{\bw_2} (X_s) - f_{\bw_1} (X_s))^2
      \\
      &&
      \quad
      +
      8  \cdot
      \frac{1}{n}
      \sum_{s=1}^n
      (Y_s-f_{\bw_2} (X_s) )^2 
      \cdot
      \sum_{k,i,j,l}
      \left(
      \frac{\partial f_{\bw_1}}{\partial w_{k,i,j}^{(l)}}(X_s)
      -
           \frac{\partial f_{\bw_2}}{\partial w_{k,i,j}^{(l)}}(X_s)
           \right)^2 
      .
    \end{eqnarray*}
    From the proof of Lemma \ref{le2} we can conclude
    \begin{eqnarray*}
      &&
      \sum_{k,i,j,l}
      \max_{s=1,\dots,n}
      \left(
      \frac{\partial f_{\bw_1}}{\partial w_{k,i,j}^{(l)}}(X_s)
      \right)^2 
      \leq
      c_{16} \cdot K_n \cdot L \cdot
(r \cdot (r+d)) \cdot r^{2L} \cdot (\gamma_n^*)^2 \cdot B_n^{2L} \cdot a^2,      
      \end{eqnarray*}
    \begin{eqnarray*}
      &&
      \frac{1}{n}
      \sum_{s=1}^n
      (f_{\bw_2} (X_s) - f_{\bw_1} (X_s))^2 
      \leq
      c_{17} \cdot K_n^2 \cdot
      (\gamma_n^*)^2 \cdot (2r+d)^{2L}
      \cdot B_n^{2L} \cdot a^2  \cdot \|\bw_1-\bw_2\|^2
      \end{eqnarray*}
    and from the proof of Lemma \ref{le3} we know
    \begin{eqnarray*}
      &&
            \frac{1}{n}
      \sum_{s=1}^n
      (Y_s-f_{\bw_2} (X_s) )^2 
      \leq
      4 \cdot K_n^2 \cdot (\gamma_n^*)^2.
    \end{eqnarray*}
    And by Lemma \ref{le2} we can conclude for any $s \in \{1, \dots, n\}$
    \[
\sum_{k,i,j,l}
      \left(
      \frac{\partial f_{\bw_1}}{\partial w_{k,i,j}^{(l)}}(X_s)
      -
           \frac{\partial f_{\bw_2}}{\partial w_{k,i,j}^{(l)}}(X_s)
           \right)^2
           \leq
           c_{18}  \cdot B_n^{4L} \cdot (\gamma_n^*)^2  \cdot
  \| \bw_1 - \bw_2 \|^2.
  \]
  Summarizing the above results we get the assertion.
    \hfill $\Box$

    By combining the above results we can show our main
    result concerning gradient descent.

    \begin{theorem}
      \label{se4th1}
      Let $A \geq 1$, $L,r \in \N$ and $K_n \in \N$.
      Define the deep neural network
  $f_\bw:\R^d \rightarrow \R$ with weight vector $\bw$
by (\ref{se2eq1})--(\ref{se2eq3}). 
Let $\beta_n,A_n, B_n \geq 1$ with $B_n \leq K_n$, and assume $X_1, \dots, X_n \in [-A,A]^d$, $|Y_i| \leq \beta_n$ $(i=1, \dots, n)$ and
\[
c_{19} \cdot K_n  \geq \beta_n \quad \mbox{and} \quad
\frac{\beta_n}{n^2} \leq c_{20}.
\]
Set
\[
F_n(\bw) = \frac{1}{n} \sum_{i=1}^n |Y_i - f_\bw(X_i)|^2.
\]
Choose some starting weight
vector $\bw^{(0)}$ which satisfies
\[
(\bw^{(0)})^{(L)}_{1,1,k}=0, \quad
|(\bw^{(0)})^{(l)}_{k,i,j}| \leq B_n \quad \mbox{and} \quad
|(\bw^{(0)})^{(0)}_{k,i,j}| \leq A_n
\]
for all $l \in \{1, \dots, L-1\}$ and all $i,j,k$,
and set
\[
\bw^{(t+1)} = \bw^{(t)} - \lambda_n \cdot \nabla_\bw F_n(\bw^{(t)})
\]
for $t=0, 1, \dots, t_n-1$. Set
\[
\lambda_n= \frac{c_{21}}{K_n^3  \cdot B_n^{2L} \cdot n}
\quad
\mbox{and}
\quad
t_n = \left \lceil
c_{22} \cdot
\frac{K_n^3  }{ \beta_n} 
\right \rceil.
\]
Let $\bw^* \in A$ where
\[
A = \left\{
\bw \quad : \quad \|\bw-\bw^{(0)}\| \leq \frac{c_{23}}{\sqrt{n} \cdot B_n^L}
\right\}
\]
and assume that (\ref{le1eq5}) holds.
Then  
\[
F_n(\bw)
\leq
F_n(\bw^*) + c_{24} \cdot \beta_n \cdot B_n^{2L} \cdot n
\cdot \|\bw^*-\bw^{(0)}\|^2
+ c_{25} \cdot \frac{\beta_n}{n}.
\]
      \end{theorem}

    \noindent
        {\bf Proof.} Set
        \[
        L_n =
        c_{26}  \cdot K_n^3  \cdot B_n^{2L} \cdot n
        \quad \mbox{and} \quad
C_n= c_{27} \cdot B_n^{2L},
\]
which implies
\[
t_n \cdot \lambda_n =
c_{28} \cdot
t_n \cdot \frac{1}{L_n}
=
c_{29} \cdot \frac{1}{\beta_n \cdot C_n \cdot n}
\quad
\mbox{and}
\quad
\beta_n \cdot \sqrt{t_n \cdot \lambda_n} \leq c_{30}. 
\]
        Because of  $(\bw^{(0)})_{1,1,k}^{(L)}=0$ we know
        $F_n(\bw^{(0)}) \leq \beta_n^2$.
        From Lemma \ref{le3} and Lemma \ref{le5}, which we apply with
        \[
        \gamma_n^* = c_{31} + c_{32} \cdot \beta_n \cdot \sqrt{t_n \cdot \lambda_n},
        \quad
        B_n= B_n +  c_{32} \cdot \beta_n \cdot \sqrt{t_n \cdot \lambda_n}
        \]
        and
        \[
        A_n= A_n +  c_{32} \cdot \beta_n \cdot \sqrt{t_n \cdot \lambda_n}
        \]
        we get
        that the assumptions of Lemma \ref{le4} are satisfied if we set
        \[
        L = c_{33} \cdot K_n^3 \cdot B_n^{2L} \cdot n.
        \]
        (In fact, $c_{34} \cdot K_n^{3/2} \cdot B_n^{2L}$ is here
        sufficient, but we use a larger value in order to get
        later that $\lambda_n \cdot D_n^2$ is small.)
        Hence
        we have
        \[
\bw^{(t)} \in A := \left\{ \bw \, : \, \|\bw-\bw^{(0)}\| \leq \frac{\sqrt{2 \beta_n^2}}{\sqrt{\beta_n \cdot C_n \cdot n}} \right\}
        \]
        for $t=1, \dots, t_n$, and
        \[
F_n( \bw^{(t_n)}) \leq \min_{t=0, \dots, t_n-1} F_n(\bw^{(t)}). 
\]
By Lemma \ref{le2} and Lemma \ref{le3} we know that the
assumptions of Lemma \ref{le1} are satisfied with
$C_n=c_{35} \cdot B_n^{2L} $ and
$D_n= c_{36} \cdot K_n^{3/2} \cdot B_n^L$.
Application of Lemma \ref{le1} with
$\delta_n=\frac{ \sqrt{2 \cdot \beta_n}}{\sqrt{\beta_n \cdot C_n \cdot n}}$
yields
\begin{eqnarray*}
  F_n( \bw^{(t_n)}) &\leq& \min_{t=0, \dots, t_n-1} F_n(\bw^{(t)}) \\
  & \leq &
    F_n(\bw^*) + \frac{\| \bw^* - \bw^{(0)} \|^2}{2 \cdot \lambda_n \cdot t_n}
  +
  12 \cdot \beta_n \cdot C_n \cdot \delta_n^2 +
  \frac{1}{2} \cdot \lambda_n \cdot D_n^2
  \\
  & \leq &
   F_n(\bw^*) + c_{38} \cdot \beta_n \cdot B_n^{2L} \cdot n
\cdot \|\bw^*-\bw^{(0)}\|^2
+ c_{39} \cdot \frac{\beta_n}{n}.
  \end{eqnarray*}
    \hfill $\Box$

    \subsection{Neural network approximation}
\label{se4sub2}
    In the sequel we construct a neural network which approximates a
    piecewise Taylor polynomial of a function $f:\R^d \rightarrow \R$.

    Assume that $f$ is $(p,C)$--smooth for some $p=q+\beta$ where
    $\beta \in (0,1]$ and $q \in \N_0$. The multivariate Taylor
      polynomial of $f$ of degree $q$ around $u \in \R^d$ is defined by
      \[
      (Tf)_{q,u}(x)
      =
      \sum_{j_1, \dots, j_d \in \N_0, \atop
        j_1 + \dots + j_d \leq q}
      \frac{
\partial^q f
      }{
\partial^{j_1} x^{(1)} \dots \partial^{j_d} x^{(d)}
      }
      (u)
      \cdot
      (x^{(1)}-u^{(1)})^{j_1} \cdot \dots \cdot (x^{(d)}-u^{(d)})^{j_d}.
      \]
Since $f$ is $(p,C)$--smooth, it is possible to show that the error
of its Taylor polynomial can be bounded by
\begin{equation}
\label{se4sub2eq1}
\left|
f(x)- (Tf)_{q,u}(x)
\right|
\leq c_{40} \cdot C \cdot \|x-u\|^p
\end{equation}
(cf., e.g., Lemma 1 in Kohler (2014)).
      For functions $f,g:\R^d \rightarrow \R$ we have
      \[
(T(f+g))_{q,u}(x)= (Tf)_{q,u}(x) + (Tg)_{q,u}(x),
      \]
      and if $g:\R^d \rightarrow \R$ is a multivariate polynomial of degree
      $q$ (or less) we have
      \[
(Tg)_{q,u}(x) = g(x) \quad (x \in \R^d).
      \]
Next we construct a piecewise Taylor polynomial.
To do this, let $A \geq 1$, $K \in \N$ and subdivide $[-A,A]^d$ into
$K^d$ many cubes of sidelength
\[
\delta=\frac{2A}{K}.
\]
Set
\[
u_k = -A + k \cdot \frac{2A}{K}
\quad (k=0, \dots, K-1).
\]
Then
\[
u_\bk = (u_{k^{(1)}}, \dots, u_{k^{(d)}})
\quad
(\bk \in I := \{0,1,\dots, K-1\}^d)
\]
denote the lower left corners of these cubes.

For $\ba, \bb \in \R^d$ we write
\[
\ba \leq \bb \quad \mbox{if} \quad a^{(l)} \leq b^{(l)} \quad
\mbox{for all } l \in \{1, \dots, d\}
\]
and
\[
\ba < \bb
 \quad \mbox{if} \quad \ba \leq \bb \quad \mbox{and} \quad \ba \neq \bb.
\]
Set
\[
  [\ba, \infty)= [a^{(1)}, \infty) \times \dots \times  [a^{(d)}, \infty)
        \quad \mbox{and} \quad
[\ba, \bb)= [a^{(1)}, b^{(1)}) \times \dots \times  [a^{(d)}, b^{(d)}).
 \]
 Our piecewise Taylor polynomial is defined by
 \[
 P(x)= \sum_{\bk \in I} P_\bk(x) \cdot 1_{
[u_\bk,\infty)
 }(x),
 \]
 where the $P_\bk$'s are recursively defined by
 \[
P_0(x) = (Tf)_{q,u_{\mathbf{0}}}(x)
 \]
 and
 \[
 P_{\bk}(x) = T \left(
f - \sum_{\bl \in I \, : \, u_\bl < u_\bk} P_\bl
 \right)_{q,u_\bk}(x).
 \]
As our next lemma shows in this way we define indeed a piecewise
Taylor polynomial.

 \begin{lemma}
   \label{le6}
   Let $\br \in I$ and $x \in [u_\br, u_\br+ \delta \cdot \mathbf{1})$.
     Then
     \[
P(x)=(Tf)_{q,u_\br}(x).
     \]
   \end{lemma}

 \noindent
     {\bf Proof.}
     The definition of $P(x)$ and $x \in [u_\br, u_\br+ \delta \cdot \mathbf{1})$
       imply
       \begin{eqnarray*}
         &&
         P(x)
         =
 \sum_{\bk \in I} P_\bk(x) \cdot 1_{
[u_\bk,\infty)
}(x)
=
\sum_{\bk \in I \, : \, u_\bk \leq u_\br} P_\bk(x)
=
P_\br(x) +
\sum_{\bk \in I \, : \, u_\bk < u_\br} P_\bk(x).
         \end{eqnarray*}
       With
       \[
 P_{\br}(x) = T \left(
f - \sum_{\bl \in I \, : \, u_\bl < u_\br} P_\bl
 \right)_{q,u_\br}(x)
 =
\left( T 
f 
 \right)_{q,u_\br}(x)
- \sum_{\bl \in I \, : \, u_\bl < u_\br} P_\bl(x)
\]
we get the assertion.
     \hfill $\Box$

     Consequently it holds
     \[
\sup_{x \in [-A,A)^d} |f(x)-P(x)| \leq c_{41} \cdot \frac{1}{K^p}
     \]
in case that $f$ is $(p,C)$--smooth (cf., (\ref{se4sub2eq1})).

In the sequel we will approximate
\[
1_{
[u_\bk,\infty)}(x)
\]
by
\[
\prod_{j=1}^d \sigma( M \cdot (x^{(j)}-u_\bk^{(j)})),
\]
where $\sigma$ is the logistic squasher and $M$ is a large
positive number. This approximation will be bad in case
that $x^{(j)}$ is close to $u_\bk^{(j)}$, and to bound the
resulting error in this case the following lemma will be useful.

\begin{lemma}
  \label{le7}
  Let $f:\R^d \rightarrow \R$ be a $(p,C)$--smooth function,
  let $\br \in I$ and let $j \in \{1, \dots, d\}$. Then we
  have for any $x \in \R^d$
with $\|u_\br-x\|_\infty \leq c_{42} \cdot \delta$:
    \[
    \left|
    \sum_{\bk \in I \, : \,
      u_\bk \leq u_\br \; and \;
      u_\bk^{(j)}=u_\br^{(j)}}
    P_\bk (x)
    \right|
    \leq
    \frac{c_{43}}{K^p}.
    \]
  \end{lemma}

\noindent
    {\bf Proof.}
Let $\be_j$ be the $j$-th unit vector in $\R^d$.
By the proof of Lemma \ref{le6}
and by (\ref{se4sub2eq1})
we have
\begin{eqnarray*}
  \left|
    \sum_{\bk \in I \, : \,
      u_\bk \leq u_\br \; and \;
      u_\bk^{(j)}=u_\br^{(j)}}
    P_\bk (x)
    \right|
    &
    =&
      \left|
    \sum_{\bk \in I \, : \,
      u_\bk \leq u_\br }
    P_\bk (x)
    -
    \sum_{\bk \in I \, : \,
      u_\bk \leq u_\br - \delta \cdot \be_j }
    P_\bk (x)
    \right|
    \\
    &
    =&
    \left|
(Tf)_{q,u_\br}(x)-(Tf)_{q,u_\br-\delta \cdot \be_j}(x)
    \right|
    \\
    &
    \leq &
       \left|
(Tf)_{q,u_\br}(x)-f(x)
       \right|
       +
          \left|
f(x)-(Tf)_{q,u_\br-\delta \cdot \be_j}(x)
\right|
\\
&
\leq &
c_{44} \cdot \|x-u_\br\|^p + c_{45} \cdot \|x-(u_{\br} - \delta \cdot \be_j)\|^p
\\
&
\leq &
\frac{c_{46}}{K^p},
  \end{eqnarray*}
where the last inequality followed from
$\|u_\br-x\|_\infty \leq c_{42} \cdot \delta= c_{42} \cdot (2A)/K$.
  \hfill $\Box$

  Next we want to approximate
  \[
 P(x)= \sum_{\bk \in I} P_\bk(x) \cdot 1_{
[u_\bk,\infty)
 }(x)
  \]
  by a neural network. Here we consider in an intermediate
  step
  \[
  \bar{P}(x)= P_0(x) + \sum_{\bk \in I, \bk \neq 0} P_\bk(x)
  \cdot
  \prod_{j=1}^d \sigma( M \cdot (x^{(j)}-u_\bk^{(j)})),
  \]
  where the indicator function is approximated by a product of neurons.

  \begin{lemma}
    \label{le8}
      Let $f:\R^d \rightarrow \R$ be a $(p,C)$--smooth function,
let $A \geq 1$, let $K \in \N$ with $K \geq e^p$, let $\sigma$
be the logistic squasher, and
      define $P(x)$ and $\bar{P}(x)$ as above.
Assume
      \[
M \geq K \cdot (\log K)^2.
\]
Then
\[
\sup_{x \in [-A,A)^d}
  | P(x) - \bar{P}(x)|
  \leq
  c_{47} \cdot \frac{1}{K^p}.
\]
    \end{lemma}

  \noindent
      {\bf Proof.}
      Let $x \in [-A,A]^d$ be arbitrary, and let $\br \in I$ be such
      that $x \in [u_\br, u_\br + \delta)$.
      Since $x \in [u_\br,u_\br+\delta) \subseteq [u_0,\infty)$
          we have
          \begin{eqnarray*}
            &&
            P(x) - \bar{P}(x)
            \\
            &&
            =
            \sum_{\bk \in I, \bk \neq 0} P_\bk(x)
            \cdot
            \left(
1_{
[u_\bk,\infty)
 }(x)            
            -
            \prod_{j=1}^d \sigma( M \cdot (x^{(j)}-u_\bk^{(j)}))
            \right)
            \\
            &&
            =
            \sum_{\bk \in I, \bk \neq 0: \atop
              u_\bk^{(j)} \geq u_\br^{(j)} + 2 \delta \mbox{ for some }
              j \in \{1, \dots, d \}
            } P_\bk(x)
            \cdot
            \left(
1_{
[u_\bk,\infty)
 }(x)            
            -
            \prod_{j=1}^d \sigma( M \cdot (x^{(j)}-u_\bk^{(j)}))
            \right)
            \\
            &&
            \quad
            +
            \sum_{\bk \in I, \bk \neq 0: \atop
              u_\bk^{(j)} \leq u_\br^{(j)} - \delta \mbox{ for all }
              j \in \{1, \dots, d \}
            } P_\bk(x)
            \cdot
            \left(
1_{
[u_\bk,\infty)
 }  (x)          
            -
            \prod_{j=1}^d \sigma( M \cdot (x^{(j)}-u_\bk^{(j)}))
            \right)
            \\
            &&
            \quad
            +
               \sum_{\bk \in I, \bk \neq 0 \; : \;
              u_\bk^{(i)} < u_\br^{(i)} + 2 \delta \mbox{ for all }
              i \in \{1, \dots, d \}, \atop
               u_\bk^{(j)} > u_\br^{(j)} - \delta \mbox{ for some }
              j \in \{1, \dots, d \}
            } P_\bk(x)
            \cdot
            \left(
1_{
[u_\bk,\infty)
 } (x)           
            -
            \prod_{j=1}^d \sigma( M \cdot (x^{(j)}-u_\bk^{(j)}))
            \right)
            \\
            &&
            =: T_{1,n}+T_{2,n}+ T_{3,n}.
            \end{eqnarray*}
          If $u_\bk^{(i)} \geq u_\br^{(i)}+2 \delta$ for some
          $i \in \{1, \dots, d\}$, then
$M \geq K \cdot (\log K)^2$ and $\delta=2A/K \geq 1/K$ imply
          \begin{eqnarray*}
            \left|
1_{
[u_\bk,\infty)
 } (x)           
            -
            \prod_{j=1}^d \sigma( M \cdot (x^{(j)}-u_\bk^{(j)}))
            \right|
            &=&
            \prod_{j=1}^d \sigma( M \cdot (x^{(j)}-u_\bk^{(j)}))
            \leq
            \sigma( M \cdot (x^{(i)}-u_\bk^{(i)}))
            \\
            &
            \leq &
            \sigma(-(\log K)^2)
            \leq
            e^{-(\log K)^2},
            \end{eqnarray*}
          which together with
          \begin{eqnarray*}
          |P_\bk(x)| &\leq& \left|
          \sum_{\bl \in I, u_{\bl} \leq u_\bk } P_\bl(x)
          -
            \sum_{\bl \in I, u_{\bl} \leq u_\bk - \delta \cdot \be_1 } P_\bl(x)
            \right|
            =\left|
 T \left(
f 
\right)_{q,u_\bk}(x)
-
 T \left(
f 
\right)_{q,u_\bk- \delta \cdot \be_1}(x)
\right|
\\
            &\leq&  c_{48}
          \end{eqnarray*}
          yields
          \[
          |T_{1,n}| \leq K^d \cdot c_{48} \cdot   e^{-(\log K)^2}
          \leq
           \frac{c_{49}}{K^p}.
           \]

           If
           $u_\bk^{(i)} \leq u_\br^{(i)}-\delta$ for all $i \in \{1, \dots, d\}$, then
           \[
           M \cdot (x^{(j)}-u_\bk^{(j)}) \geq M \cdot \delta \geq (\log K)^2
           \quad \mbox{for all } j \in \{1, \dots,d\},
           \]
           hence
           \begin{eqnarray*}
             &&
             \left|
1_{
[u_\bk,\infty)
 }            (x)
            -
            \prod_{j=1}^d \sigma( M \cdot (x^{(j)}-u_\bk^{(j)}))
            \right|
            \\
            &&
            =
            1 - \prod_{j=1}^d \sigma( M \cdot (x^{(j)}-u_\bk^{(j)}))
            \\
            &&
            =
            1 - \prod_{j=1}^d \frac{1}{1 + e^{-M \cdot (x^{(j)}-u_\bk^{(j)})}}
            \\
            &&
           \leq
            1 - \prod_{j=1}^d \frac{1}{1 + e^{-(\log K)^2}}
            \\
            &&
            =
            \sum_{l=1}^d
            \left(
            \prod_{j=1}^{l-1} \frac{1}{1 + e^{-(\log K)^2}}
            -
            \prod_{j=1}^l \frac{1}{1 + e^{-(\log K)^2}}
            \right)
            \\
            &&
            \leq
            d \cdot
            \left(
1-\frac{1}{1 + e^{-(\log K)^2}}
\right)
\\
&&
\leq
d \cdot  e^{-(\log K)^2},
             \end{eqnarray*}
           which implies
           \[
           |T_{2,n}| \leq
K^d \cdot c_{48} \cdot  d \cdot e^{-(\log K)^2}
          \leq
           \frac{c_{50}}{K^p}.           
           \]
           So it remains to bound $|T_{3,n}|$.

           $T_{3,n}$ is a sum of less than
           \[
3^d
           \]
           terms of the form
           \[
               \sum_{\bk \in I, \bk \neq 0 \; : \;
              u_\bk^{(i)} \leq u_\br^{(i)} - \delta \mbox{ for all }
              i \in \{1, \dots, d \} \setminus \{j_1, \dots, j_s \},\atop
               u_\bk^{(j_t)} = u_\br^{(j_t)} + l_t \cdot \delta \mbox{ for all }
              t \in \{1, \dots, s \}
            } P_\bk(x)
            \cdot
            \left(
1_{
[u_\bk,\infty)
}
(x)
            -
            \prod_{j=1}^d \sigma( M \cdot (x^{(j)}-u_\bk^{(j)}))
            \right)
,
           \]
           where
           $s \in \{1, \dots, d\}$, $1 \leq j_1 < j_2 < \dots < j_s \leq d$,
           $l_1, \dots, l_s \in \{0,1\}$.
           The absolute value of the difference of this term and the term
           \[
     \sum_{\bk \in I, \bk \neq 0 \; : \;
              u_\bk^{(i)} \leq u_\br^{(i)} - \delta \mbox{ for all }
              i \in \{1, \dots, d \} \setminus \{j_1, \dots, j_s \},
              \atop
               u_\bk^{(j_t)} = u_\br^{(j_t)} + l_t \cdot \delta \mbox{ for all }
              t \in \{1, \dots, s \}
            } P_\bk(x)
            \cdot
            \left(
1_{
[u_\bk,\infty)
}
(x)
            -
            \prod_{t=1}^s \sigma( M \cdot (x^{(j_t)}-u_\bk^{(j_t)}))
            \right)
           \]
           is because of
           \begin{eqnarray*}
             &&
           \left|
            \prod_{j=1}^d \sigma( M \cdot (x^{(j)}-u_\bk^{(j)}))
                       -
            \prod_{t=1}^s \sigma( M \cdot (x^{(j_t)}-u_\bk^{(j_t)}))
            \right|
            \\
            &&
            \leq
            \left|
            1 -
            \prod_{j \in \{1, \dots, d\} \setminus
            \{j_1, \dots, j_s\}} \sigma( M \cdot (x^{(j)}-u_\bk^{(j)}))
            \right|
            \\
            &&
           \leq
           d \cdot e^{- (\log K)^2}
           \end{eqnarray*}
           (which follows as above) bounded from above by $c_{51}/K^p$.

           Hence it suffices to show that
           \begin{equation}
             \label{ple8eq1}
    \sum_{\bk \in I, \bk \neq 0 \; : \;
              u_\bk^{(i)} \leq u_\br^{(i)} - \delta \mbox{ for all }
              i \in \{1, \dots, d \} \setminus \{j_1, \dots, j_s \},
              \atop
               u_\bk^{(j_t)} = u_\br^{(j_t)} + l_t \cdot \delta \mbox{ for all }
              t \in \{1, \dots, s \}
            } P_\bk(x)
            \cdot
1_{
[u_\bk,\infty)
 }(x)            
\end{equation}
and
\[
     \sum_{\bk \in I, \bk \neq 0 \; : \;
              u_\bk^{(i)} \leq u_\br^{(i)} - \delta \mbox{ for all }
              i \in \{1, \dots, d \} \setminus \{j_1, \dots, j_s \},
              \atop
               u_\bk^{(j_t)} = u_\br^{(j_t)} + l_t \cdot \delta \mbox{ for all }
              t \in \{1, \dots, s \}
            } P_\bk(x)
            \cdot
            \prod_{t=1}^s \sigma( M \cdot (x^{(j_t)}-u_\bk^{(j_t)}))
\]
are bounded in absolute value by $c_{52}/K^p$. Since
\begin{eqnarray*}
  &&
  \sum_{\bk \in I, \bk \neq 0 \; : \;
              u_\bk^{(i)} \leq u_\br^{(i)} - \delta \mbox{ for all }
              i \in \{1, \dots, d \} \setminus \{j_1, \dots, j_s \},
              \atop
               u_\bk^{(j_t)} = u_\br^{(j_t)} + l_t \cdot \delta \mbox{ for all }
              t \in \{1, \dots, s \}
            } P_\bk(x)
            \cdot
            \prod_{t=1}^s \sigma( M \cdot (x^{(j_t)}-u_\bk^{(j_t)}))
            \\
            &&
            =
            \prod_{t=1}^s \sigma( M \cdot (x^{(j_t)}- u_\br^{(j_t)} - l_t \cdot \delta ))
            \cdot
            \sum_{\bk \in I, \bk \neq 0 \; : \;
              u_\bk^{(i)} \leq u_\br^{(i)} - \delta \mbox{ for all }
              i \in \{1, \dots, d \} \setminus \{j_1, \dots, j_s \},
              \atop
               u_\bk^{(j_t)} = u_\br^{(j_t)} + l_t \cdot \delta \mbox{ for all }
              t \in \{1, \dots, s \}
            } P_\bk(x)
  \end{eqnarray*}
for this it suffices to show that terms of the form
\[
 \sum_{\bk \in I, \bk \neq 0 \; : \;
              u_\bk^{(i)} \leq u_\br^{(i)} \mbox{ for all }
              i \in \{1, \dots, d \} \setminus \{j_1, \dots, j_s \}, \atop
               u_\bk^{(j_t)} = u_\br^{(j_t)}  \mbox{ for all }
              t \in \{1, \dots, s \}
            } P_\bk(x),
\]
where $u_\br$ satisfies $\|x-u_{\br}\|_\infty \leq 2 \delta$, are
bounded in absolute value by 
$c_{53}/K^p$, which we do in the sequel.
(Here we have used that in
(\ref{ple8eq1})
w.l.o.g.
$1_{[u_\bk,\infty)}(x)=1$
  holds because otherwise (\ref{ple8eq1}) is zero.)

Since
\begin{eqnarray*}
  &&
 \sum_{\bk \in I, \bk \neq 0 \; : \;
              u_\bk^{(i)} \leq u_\br^{(i)} \mbox{ for all }
              i \in \{1, \dots, d \} \setminus \{j_1, \dots, j_s \}\atop
               u_\bk^{(j_t)} = u_\br^{(j_t)}  \mbox{ for all }
              t \in \{1, \dots, s \}
 } P_\bk(x)
 \\
 &&
 =
  \sum_{\bk \in I, \bk \neq 0 \; : \;
              u_\bk \leq u_\br, \atop
               u_\bk^{(j_t)} = u_\br^{(j_t)}  \mbox{ for all }
              t \in \{2, \dots, s \}
 } P_\bk(x)
  -
   \sum_{\bk \in I, \bk \neq 0 \; : \;
              u_\bk \leq u_\br- \be_{j_1} \cdot \delta, \atop
               u_\bk^{(j_t)} = u_\br^{(j_t)}  \mbox{ for all }
              t \in \{2, \dots, s \}
 } P_\bk(x)
  \end{eqnarray*}
we see that the term above is equal to a sum of at most $2^{s-1}$ terms
of the form
\[
    \sum_{\bk \in I \, : \,
      u_\bk \leq u_\br  \mbox{ and } 
      u_\bk^{(j)}=u_\br^{(j)}}
    P_\bk (x),
    \]
    where $\|u_\br-x\|_\infty \leq 3 \cdot \delta$.
From this the assertion follows by an application of Lemma \ref{le7}.
\hspace*{3cm} \hfill $\Box$

Next we want to approximate
  \[
  \bar{P}(x)= P_0(x) + \sum_{\bk \in I, \bk \neq 0} P_\bk(x)
  \cdot
  \prod_{j=1}^d \sigma( M \cdot (x^{(j)}-u_\bk^{(j)}))
  \]
  by a neural network. In order to do this, we need to represent monomials
  by neural networks and need to be able to multiply real numbers
  by using neural networks. The starting point for both is the following
  lemma, which is a modification of Theorem 2 in Scarselli and Tsoi (1998).

  \begin{lemma}
    \label{le9}
    Let $\sigma$ be the logistic squasher, let $k \in \N$, let $t_\sigma \in \R$
    be such that $\sigma^{(k)}(t_\sigma) \neq 0$.
    Then for any $N \in \N$ with $N>k$ there exist
    \[
\alpha_j, \beta_j \in \R \quad (j=0, \dots, N-1)
\]
such that
\[
f_{net,x^k}(x)
=
\frac{k!}{
\sigma^{(k)}(t_\sigma)
}
\cdot
\sum_{j=0}^{N-1}
\alpha_j \cdot \sigma( \beta_j \cdot x + t_\sigma)
\]
satisfies for all $A>0$ and all $x \in [-A,A]$:
\[
\left|
f_{net,x^k}(x)
-x^k
\right|
\leq c_{54} \cdot A^N
\]
for some $c_{54}=c_{54}(N,k,\sigma^{(k)}(t_\sigma),
\| \sigma^{(N)} \|_\infty, \alpha_0, \dots, \alpha_{N-1},\beta_0, \dots,  \beta_{N-1}) \geq 0$.
    \end{lemma}

  \noindent
      {\bf Proof.} Let $\beta_j \in \R$ $(j=0, \dots, N-1)$ be
      pairwise distinct. Then the vectors
      \[
      \bv_l = (\beta_0^l, \dots, \beta_{N-1}^l)^T
      \quad
      (l=0, \dots, N-1)
      \]
      are linearly independent since
      \[
\sum_{l=0}^{N-1} \alpha_l \cdot \bv_l=0
      \]
      implies that the polynomial
      \[
p(x)= \sum_{l=0}^{N-1} \alpha_l \cdot x^l
      \]
      of degree $N-1$ has the $N$ zero points $\beta_0$, \dots, $\beta_{N-1}$,
      which is possible only in case $\alpha_0=\dots=\alpha_{N-1}=0$.
      Hence we can choose $\alpha_0, \dots, \alpha_{N-1} \in \R$ such
      that
      \[
\alpha_0 \cdot \bv_0 + \dots + \alpha_{N-1} \cdot \bv_{N-1}
      \]
      is equal to the $k$-th unit vector in $\R^N$, which implies
      \begin{equation}
        \label{ple9eq1}
        \sum_{j=0}^{N-1}
        \alpha_j \cdot \beta_j^l
        =
        \begin{cases}
          1, & \mbox{if } l=k \\
          0, & \mbox{if } l \in \{0, \dots, N-1\} \setminus \{k\}.
          \end{cases}
      \end{equation}
      Using these values for the $\alpha_j$ and $\beta_j$, a Taylor
      expansion of
      \[
x \mapsto \sigma(\beta_j \cdot x + t_\sigma)
      \]
      around $t_\sigma$ of order $N-1$ implies
      \begin{eqnarray*}
        f_{net,x^k}(x) &=&
        \frac{k!}{
\sigma^{(k)}(t_\sigma)
}
\cdot
\sum_{j=0}^{N-1}
\alpha_j \cdot
\left(
\sum_{l=0}^{N-1}
\frac{\sigma^{(l)}(t_\sigma)}{l!}
\cdot
(\beta_j \cdot x)^l
+
\frac{\sigma^{(N)}(\xi_j)}{l!}
\cdot
(\beta_j \cdot x)^N
\right)
\\
&=&       \frac{k!}{
\sigma^{(k)}(t_\sigma)
}
\cdot
\sum_{l=0}^{N-1}
\left(
\sum_{j=0}^{N-1} \alpha_j \cdot \beta_j^l
\right)
\cdot
\frac{\sigma^{(l)}(t_\sigma)}{l!} \cdot x^l
\\
&&
\quad
+
 \frac{k!}{
\sigma^{(k)}(t_\sigma)
}
 \cdot
 \sum_{j=0}^{N-1}
 \alpha_j \cdot
 \frac{\sigma^{(N)}(\xi_j)}{l!}
 \cdot 
 \beta_j^N \cdot x^N
 \\
 &=&
 x^k +  \frac{k!}{
\sigma^{(k)}(t_\sigma)
}
 \cdot
 \sum_{j=0}^{N-1}
 \alpha_j \cdot
 \frac{\sigma^{(N)}(\xi_j)}{l!}
 \cdot 
 \beta_j^N \cdot x^N,
      \end{eqnarray*}
      where the last equality follows from (\ref{ple9eq1}).
      Hence
      \[
\left|
f_{net,x^k}(x)
-x^k
\right|
\leq
\left|
 \frac{k!}{
\sigma^{(k)}(t_\sigma)
}
 \cdot
 \sum_{j=0}^{N-1}
 \alpha_j \cdot 
 \frac{\sigma^{(N)}(\xi_j)}{l!}
 \cdot 
 \beta_j^N \cdot x^N
 \right|
 \leq c_{55} \cdot |x|^N
 \leq c_{55} \cdot A^N.
      \]
      \hfill $\Box$

      Our next lemma uses Lemma \ref{le9} in order
      to construct a neural network which can multiply two numbers.

      \begin{lemma}
        \label{le10}
    Let $A>0$, let $N \in \N$ with $N>2$ and let $f_{net,x^2}$ be the neural
    network from Lemma \ref{le9} (which has one hidden layer with $N$
neurons). Then
    \[
    f_{mult}(x,y)= \frac{1}{4} \cdot \left(
f_{net,x^2}(x+y) - f_{net,x^2}(x-y)
    \right)
    \]
    satisfies for all $x,y \in [-A,A]$:
\[
\left|
f_{mult}(x,y) - x \cdot y
\right|
\leq c_{56} \cdot A^N
\]
for some $c_{56} > 0$ (which depends on the constant $c_{54}$ in Lemma \ref{le9}
and on $N$).
        \end{lemma}

      \noindent
          {\bf Proof.} By Lemma \ref{le9} we get
          \begin{eqnarray*}
            &&
            \left|
f_{mult}(x,y) - x \cdot y
\right|
\\
&&
=
\left|
\frac{1}{4} \cdot \left(
f_{net,x^2}(x+y) - f_{net,x^2}(x-y)
\right)
-
\frac{1}{4} \cdot \left(
(x+y)^2 - (x-y)^2
    \right)
    \right|
    \\
    &&
    \leq
    \frac{1}{4} \cdot | f_{net,x^2}(x+y) - (x+y)^2|
    +
    \frac{1}{4} \cdot | f_{net,x^2}(x-y) - (x-y)^2|
    \\
    &&
    \leq
    \frac{1}{4} \cdot
    c_{54} \cdot (2A)^N
    +
      \frac{1}{4} \cdot
      c_{54} \cdot (2A)^N
      \leq
      c_{56} \cdot A^N.
          \end{eqnarray*}
          \hfill $\Box$

          Next we extend the multiplication network from the previous lemma
          in such a way that it can multiply a finite number of real
          values simultaneously.

          \begin{lemma}
            \label{le11}
            Let $\sigma(x)=1/(1+\exp(-x))$, let $0<A \leq 1$, let $N \in \N$ with $N>2$ and
            let $d \in \N$.
            Assume
            \begin{equation}
              \label{le11eq1}
              c_{56} \cdot 4^{d \cdot N} \cdot A^{N-1} \leq 1,
            \end{equation}
            where $c_{56}$ is the constant from Lemma \ref{le10}.
            Then there exists a neural network
            \[
f_{mult,d}
\]
with at most $\lceil \log_2 d \rceil$ many layers, at most
$2 \cdot N \cdot d$
many neurons and activation function $\sigma$,
where all the weights are bounded in absolute value by some constant,
such
that for all $x_1, \dots, x_d \in [-A,A]$ it holds:
\[
| f_{mult,d}(x_1, \dots, x_d) - \prod_{j=1}^d x_j| \leq c_{57} \cdot A^N,
\]
where $c_{57} \geq 0$.
            \end{lemma}

          \noindent
              {\bf Proof.} The proof is a modification of the proof of Lemma 7
              in Kohler and Langer (2021).

              We set $q=\lceil \log_2(d)\rceil$. The feedforward neural network $f_{mult, d}$ with $L=q$ hidden layers and $r=
2 \cdot N \cdot d$
neurons in each layer is constructed as follows: Set 
\begin{equation}
\label{neq1}
(z_1, \dots, z_{2^q})=
  \left(x^{(1)}, x^{(2)}, \dots, x^{(d)}, \underbrace{1, \dots,1}_{2^q-d \; times} \right).
\end{equation}
In the construction of our network we will use the network $f_{mult}$ of Lemma \ref{le10},  
which satisfies
\begin{equation}
  \label{ple11eq1}
  |f_{mult}(x,y) - x \cdot y| \leq c_{56} \cdot 4^{d \cdot N} \cdot A^N
\end{equation}
for $x,y \in [-4^{d} A,4^{d} A]$.
In the first layer we compute
\[
f_{mult}(z_1,z_2), 
f_{mult}(z_3,z_4), 
\dots,
f_{mult}(z_{2^q-1},z_{2^q}), 
\]
which can be done by one layer of
$2 \cdot N \cdot 2^{q-1} \leq 2 \cdot N \cdot d$
neurons. 
As a result of the first layer we get a vector of outputs
which has length $2^{q-1}$. Next we pair these outputs and apply $f_{mult}$ again. This procedure is continued until there is only one output left.
Therefore we need $L =q$ hidden layers and
at most
$2 \cdot N \cdot d$
neurons in each layer. 

By  (\ref{le11eq1}) and (\ref{ple11eq1}) 
we get for any $l \in \{1,\dots,d\}$ and any
$z_1,z_2 \in [-(4^l-1) \cdot A,(4^l-1) \cdot A]$
\[
|\hat{f}_{mult}(z_1,z_2)| \leq
|z_1 \cdot z_2| + |\hat{f}_{mult}(z_1,z_2)-z_1 \cdot z_2|
\leq
(4^l-1)^2 A^{2l} + c_{56} \cdot 4^{d \cdot N} \cdot A^N
\leq
(4^{2l}-1) \cdot A.
\]
From this we get successively that all outputs
of 
layer $l \in \{1,\dots,q-1\}$
  are contained in the interval
$[-(4^{2^l}-1) \cdot A,(4^{2^l}-1) \cdot A]$, hence in particular they
  are contained in the interval
$[-4^{d} A,4^{d} A]$
where inequality  (\ref{ple11eq1}) does hold.

  Define $\hat{f}_{2^q}$ recursively by
\[
\hat{f}_{2^q}(z_1,\dots,z_{2^q})
=
\hat{f}_{mult}(\hat{f}_{2^{q-1}}(z_1,\dots,z_{2^{q-1}}),\hat{f}_{2^{q-1}}(z_{2^{q-1}+1},\dots,z_{2^q}))
\]
and
\[
\hat{f}_2(z_1,z_{2})= \hat{f}_{mult}(z_1,z_{2}),
\]
and set
\[
\Delta_l=\sup_{z_1,\dots,z_{2^l}
  \in [-A,A]}
|\hat{f}_{2^l}(z_1,\dots,z_{2^l})-  \prod_{i=1}^{2^l} z_i|.
\]
Then
\[
|\hat{f}_{mult, d}(x_1, \dots, x_d)-\prod_{i=1}^d x_i|
\leq
\Delta_q
\]
and from
\[
\Delta_1 \leq c_{56} \cdot 4^{d \cdot N} \cdot A^{N}
\]
(which follows from (\ref{ple11eq1})) and
\begin{eqnarray*}
  \Delta_q
  & \leq &
  \sup_{z_1,\dots,z_{2^q}
  \in [-A,A]}
  |\hat{f}_{mult}(\hat{f}_{2^{q-1}}(z_1,\dots,z_{2^{q-1}}),\hat{f}_{2^{q-1}}(z_{2^{q-1}+1},\dots,z_{2^q}))
  \\
  &&
  \hspace*{4cm}
    -
    \hat{f}_{2^{q-1}}(z_1,\dots,z_{2^{q-1}}) \cdot \hat{f}_{2^{q-1}}(z_{2^{q-1}+1},\dots,z_{2^q})|
      \\
      &&
      +
  \sup_{z_1,\dots,z_{2^q}
  \in [-A,A]}
      \left|\hat{f}_{2^{q-1}}(z_1,\dots,z_{2^{q-1}}) \cdot \hat{f}_{2^{q-1}}(z_{2^{q-1}+1},\dots,z_{2^q})\right.
      \\
      &&
      \hspace*{4cm}
      -
        \left.\left( \prod_{i=1}^{2^{q-1}} z_i \right)
\cdot \hat{f}_{2^{q-1}}(z_{2^{q-1}+1},\dots,z_{2^q})\right|
      \\
      &&
      +
  \sup_{z_1,\dots,z_{2^q}
  \in [-A,A]}
      \left|
        \left( \prod_{i=1}^{2^{q-1}} z_i \right)
        \cdot \hat{f}_{2^{q-1}}(z_{2^{q-1}+1},\dots,z_{2^q})
          \right.
          \\
          && \left. 
          \hspace*{6cm}
          -
          \left( \prod_{i=1}^{2^{q-1}} z_i \right)
          \cdot
          \prod_{i=2^{q-1}+1}^{2^{q}} z_i
          \right|
          \\
          &
          \leq &
          c_{56} \cdot 4^{d \cdot N} \cdot A^N
          + 4^{2^{q-1}} \cdot A \cdot
          \Delta_{q-1}
          +  A^{2^{q-1}} \cdot
          \Delta_{q-1}
          \\
          &\leq&
           c_{56} \cdot 4^{d \cdot N} \cdot A^N
          + 2 \cdot 4^{2^{q-1}}  \cdot
          \Delta_{q-1}
  \end{eqnarray*}
(where  the second inequality follows from
(\ref{ple11eq1})
and the fact that all outputs of 
layer $l \in \{1,\dots,q-1\}$
  are contained in the interval
$[-4^{2^l} A,4^{2^l} A]$)
we get
for $x \in [-A,A]^d$
\begin{eqnarray*}
  |\hat{f}_{mult, d}(\bold{x}) - \prod_{i=1}^d x^{(i)}|
  & \leq & \Delta_q \nonumber \\
  &
  \leq &
             c_{56} \cdot 4^{d \cdot N} \cdot A^N \cdot
  4^{1+2+\dots+2^{q-1}}  \cdot
  \left(
1 + 2 + \dots + 2^{q-1}
  \right)
  \nonumber \\
  &
  \leq &
           c_{56} \cdot 4^{d \cdot N} \cdot A^N \cdot 4^{2d+1} \cdot d
  \nonumber \\
  &
  = &
             c_{56} \cdot 4^{d \cdot N+2d+1} \cdot d \cdot A^N 
  ,
\end{eqnarray*}
where the last inequality was implied by
\[
1+2+ \dots + 2^{q-1} = 2^q \leq 2 \cdot d.
\]
\hfill $\Box$

We are now ready to formulate and prove our main result about
the approximation of $(p,C)$--smooth function by deep neural
networks with bounded weights.

\begin{theorem}
  \label{se4th2}
Let $d \in \N$,  $p=q+\beta$ where
$\beta \in (0,1]$ and $q \in \N_0$, $C>0$,
$A \geq 1$
and
$A_n, B_n, \gamma_n^* \geq 1$.
For $L,r,K \in \N$
let $\F$ be the set of all networks $f_{\bw}$ defined by
(\ref{se2eq1})--(\ref{se2eq3}) with $K_n$ replaced by $r$, where
the weight vector satisfies
\[
|w_{i,j}^{(0)}| \leq A_n, \quad
|w_{i,j}^{(l)}| \leq B_n \quad \mbox{and} \quad
|w_{i,j}^{(L)}| \leq \gamma_n^*
\]
for all $l \in \{1, \dots, L-1\}$ and all $i,j$, and set
\[
\HH = \left\{ \sum_{k=1}^{K^d} f_k \quad : \quad f_k \in \F \quad (k=1, \dots, K)
\right\}.
\]
Let $L,r \in \N$ with
\[
L \geq \lceil \log_2(q+d) \rceil
\quad
\mbox{and}
\quad
r \geq 2 \cdot (2p+d) \cdot (q+d),
\]
and set
\[
 A_n =A \cdot K \cdot \log K, \quad B_n=c_{58}
\quad \mbox{and} \quad 
\gamma_n^*=c_{59} \cdot K^{q+d}.
\]
Assume $K \geq c_{60}$ for $c_{60}>0$ sufficiently large.
Then there exists for any $(p,C)$--smooth $f:\R^d \rightarrow \R$
a neural network $h \in \HH$ such that
\[
\sup_{x \in [-A,A)^d} |f(x)-h(x)|
\leq
\frac{c_{61}}{K^{p}}.
\]
\end{theorem}

\noindent
    {\bf Proof.}
    Define $P(x)$ and $\bar{P}(x)$ as above with $M=K \cdot (\log K)^2$.
    Then
     \[
\sup_{x \in [-A,A)^d} |f(x)-P(x)| \leq c_{62} \cdot \frac{1}{K^p}
     \]
     and
\[
\sup_{x \in [-A,A)^d}
  | P(x) - \bar{P}(x)|
  \leq
  c_{63} \cdot \frac{1}{K^p}
\]
(cf., Lemma \ref{le8}) imply that it suffices to show
that there exists $h \in \HH$ such that
\[
\sup_{x \in [-A,A)^d}
  | h(x) - \bar{P}(x)|
  \leq
  c_{64} \cdot \frac{1}{K^p}.
\]
Since $\bar{P}(x)$ is a sum of $P_0(x)$ and $(K^d-1)$ terms of the form
\[
P_\bk(x)
  \cdot
  \prod_{j=1}^d \sigma( M \cdot (x^{(j)}-u_\bk^{(j)})),
\]
where each $P_k(x)$ is a multivariate polynomial polynomial of
degree $q$ with bounded coefficients, if suffices to show
that for all $i_1, \dots, i_q \in \{0,\dots,d\}$, all
$u \in [-A,A]^d$ and
$z_0=1$, $z_j=x^{(j)}-u^{(j)}$
$(j=1, \dots, d)$
there exists $f_1, f_2 \in \F$ such that
\begin{equation}
  \label{pse4th2eq1}
\sup_{x \in [-A,A]^d}
| \prod_{s=1}^q z_{i_s} 
-
f_1(x)|
\leq
\frac{c_{65}}{K^{p+d}}
\end{equation}
and
\begin{equation}
  \label{pse4th2eq2}
\sup_{x \in [-A,A]^d}
| \prod_{s=1}^q z_{i_s} \cdot
\prod_{j=1}^d \sigma( M \cdot (x^{(j)}-u_\bk^{(j)}))
-
f_2(x)|
\leq
\frac{c_{66}}{K^{p+d}}.
\end{equation}
Let $f_{id}=f_{net,x}$ be the network of Lemma \ref{le9} which satisfies
\begin{equation}
\label{pse4th2eq3}
|f_{id}(x)-x| \leq
\frac{c_{67}}{K^{2p+2d}}
\end{equation} 
for all
$x \in [-c_{68}/K,c_{68}/K]$ (so we use $N=\lceil 2p+2d \rceil$).
Set
\[
f_{id}^{(1)}=f_{id} \quad \mbox{and} \quad f_{id}^{(l+1)}=f_{id}^{(l)} \circ f_{id}
\]
for $l \in \N$.  Because of (\ref{pse4th2eq3}) and
\begin{eqnarray*}
| f_{id}^{(l+1)}(x)-x|
& \leq &
| f_{id}^{(l+1)}(x)-f_{id}^{(l)}(x)|
+
| f_{id}^{(l)}(x)-x|
\end{eqnarray*}
an easy induction shows
\begin{equation}
\label{pse4th2eq4}
|f_{id}^{(l)}(x)-x| \leq
\frac{c_{69,l}}{K^{2p+2d}}
\end{equation} 
for all $x \in [-c_{70}/K,c_{70}/K]$.

Furthermore, let $f_{mult,q}$ and $f_{mult,q+d}$ be the networks
from Lemma \ref{le11} which satisfy
\[
| f_{mult,q}(x_1, \dots, x_q) - \prod_{j=1}^q x_j| \leq 
\frac{c_{71}}{K^{2p+d}},
\]
for all $x_1, \dots, x_q \in [-c_{72}/K,c_{72}/K]$ and
\[
| f_{mult,q+d}(x_1, \dots, x_{q+d}) - \prod_{j=1}^{q+d} x_j| \leq 
\frac{c_{73}}{K^{2p+2d}},
\]
for all $x_1, \dots, x_{q+d} \in [-c_{74}/K,c_{74}/K]$.

Then we define
\[
f_1(x)= K^{q} \cdot f_{id}^{(L-\lceil \log_2q \rceil)} ( f_{mult,q}(
z_{i_1}/K, \dots, z_{i_q}/K
))
\]
and
\begin{eqnarray*}
  f_2(x)&=&K^{q+d} \cdot
  f_{id}^{(L-\lceil \log_2(q+d) \rceil)} \Bigg(
f_{mult,q+d} \Bigg(f_{id}(z_{i_1}/K), \dots, f_{id}(z_{i_q}/K),
\\
&&
\hspace*{1cm}
\frac{1}{K} \cdot \sigma( M \cdot (x^{(1)}-u^{(1)})),
\dots,
\frac{1}{K} \cdot \sigma( M \cdot (x^{(d)}-u^{(d)})) \Bigg) \Bigg).
\end{eqnarray*}
Then $f_1$ and $f_2$ are both contained in $\F$.
Using (\ref{pse4th2eq1}) and (\ref{pse4th2eq4}) we get
\begin{eqnarray*}
  &&
| \prod_{s=1}^q z_{i_s} 
-
f_1(x)|
\\
&&
\leq
K^q \cdot
\left|
\prod_{s=1}^q z_{i_s}/K^q 
-
f_{id}^{(L-\lceil \log_2q \rceil)}
(
f_{mult,q}(
z_{i_1}/K, \dots, z_{i_q}/K
))
\right|
\\
&&
\leq
K^q \cdot
\left(
\left|
\prod_{s=1}^q z_{i_s}/K^q 
-
 f_{mult,q}(
z_{i_1}/K, \dots, z_{i_q}/K
)
\right|
+
\frac{c_{75}}{K^{2p+2d}}
\right)
\\
&&
\leq
K^q \cdot
\left(
\frac{c_{76}}{K^{2p+2d}}
+
\frac{c_{75}}{K^{2p+2d}}
\right)
\leq
\frac{c_{77}}{K^{p+d}}
\end{eqnarray*} 
and 
\begin{eqnarray*}
  &&
| \prod_{s=1}^q z_{i_s} \cdot
\prod_{j=1}^d \sigma( M \cdot (x^{(j)}-u_\bk^{(j)}))
-
f_2(x)|
\\
&&
=
K^{q+d} \cdot
\Bigg|
\prod_{s=1}^q z_{i_s}/K \cdot
\prod_{j=1}^d \frac{1}{K} \cdot \sigma( M \cdot (x^{(j)}-u_\bk^{(j)}))
\\
&&
\quad
-
f_{id}^{(L- \lceil \log_2(q+d) \rceil)}
\Bigg(
f_{mult,q+d} \Bigg(f_{id}(z_{i_1}/K), \\
&&
\hspace*{1cm}
\dots, f_{id}(z_{i_q}/K),
\sigma( M \cdot (x^{(1)}-u^{(1)}))/K,
\dots,
\sigma( M \cdot (x^{(d)}-u^{(d)})/K) \Bigg) \Bigg)
\Bigg|
\\
&&
\leq
K^{q+d} \cdot
\Bigg( \frac{c_{78}}{K^{2p+2d}}
+
\Bigg|
\prod_{s=1}^q f_{id}(z_{i_s}/K) \cdot
\prod_{j=1}^d \frac{1}{K} \cdot \sigma( M \cdot (x^{(j)}-u_\bk^{(j)}))
\\
&&
\hspace*{2cm}
-
f_{mult,q+d} \Bigg(f_{id}(z_{i_1}/K), \dots, f_{id}(z_{i_q}/K),
\sigma( M \cdot (x^{(1)}-u^{(1)}))/K,
\\
&&
\hspace*{3cm}
\dots,
\sigma( M \cdot (x^{(d)}-u^{(d)})/K) \Bigg)
\Bigg|
\Bigg)
\\
&&
\leq
K^{q+d} \cdot
\left( \frac{c_{78}}{K^{2p+2d}}
+\frac{c_{79}}{K^{2p+2d}} \right)
\leq
\frac{c_{80}}{K^{p+d}},
\end{eqnarray*}
which implies the assertion. \hfill $\Box$

\subsection{Neural network generalization}
\label{se4sub3}
In order to control the generalization error
of our over-parameterized spcaes of deep neural networks
we use the following metric entropy bound.

\begin{lemma}
  \label{le12} 
  Let $\alpha, \beta \geq 1$ and let $A,B,C \geq 1$.
  Let $\sigma:\R \rightarrow \R$ be $k$-times differentiable
  such that all derivatives up to order $k$ are bounded on $\R$.
  Let $\F$
  be the set of all functions $f_{\bw}$ defined by
  (\ref{se2eq1})--(\ref{se2eq3}) where the weight vector $\bw$
  satisfies
  \begin{equation}
    \label{le12eq1}
    \sum_{j=1}^{K_n} |w_{1,1,j}^{(L)}| \leq C,
    \end{equation}
  \begin{equation}
    \label{le12eq2}
    |w_{k,i,j}^{(l)}| \leq B \quad (k \in \{1, \dots, K_n\},
    i,j \in \{1, \dots, r\}, l \in \{1, \dots, L-1\})
    \end{equation}
and
  \begin{equation}
    \label{le12eq3}
    |w_{k,i,j}^{(0)}| \leq A \quad (k \in \{1, \dots, K_n\},
    i \in \{1, \dots, r\}, j \in \{1, \dots,d\}).
  \end{equation}
  Then we have for any $1 \leq p < \infty$, $0 < \epsilon < 1$ and
  $x_1^n \in \Rd$
  \begin{eqnarray*}
&&\Nu_p \left(
\epsilon, \{ T_\beta f \cdot 1_{[-\alpha,\alpha]^d} \, : \, f \in \F \}, x_1^n
\right)
\\
&&
\leq \left(c_{81}\cdot \frac{\beta^p} {\epsilon^p}\right)^{c_{82}\cdot \alpha^d \cdot B^{(L-1)\cdot d} \cdot A^d \cdot \left(\frac{C}{\epsilon}\right)^{d/k}+ c_{83}
}.
\\
  \end{eqnarray*}
  
  \end{lemma}

\noindent
    {\bf Proof.} This result follows from Lemma 4 in Drews and Kohler (2024). For the sake of completeness
    we repeat the proof below.
    
    In the {\it first step} of the proof we show
    for any $f_{\bw} \in \F$, any $x \in \Rd$ and any $s_1, \dots, s_k \in \{1, \dots, d\}$
    \begin{equation}
      \label{ple12eq1}
      \left|
      \frac{\partial^k f_{\bw}}{\partial x^{(s_1)} \dots \partial x^{(s_k)}} (x)
      \right| \leq c_{84} \cdot C \cdot B^{(L-1) \cdot k} \cdot A^k =: c
      .
      \end{equation}
    The definition of $f_{\bw}$ implies
    \[
    \frac{\partial^k f_{\bw}}{\partial x^{(s_1)} \dots \partial x^{(s_k)}} (x) = \sum_{j=1}^{K_n} w_{1,1,j}^{(L)} \cdot
    \frac{\partial^k
f_{j,1}^{(L)}(x)
    }{\partial x^{(s_1)} \dots \partial x^{(s_k)}} (x),
    \]
    hence (\ref{ple12eq1}) is implied by
    \begin{equation}
      \label{ple12eq1b}
      \left|
      \frac{\partial^k f_{j,1}^{(L)}}{\partial x^{(s_1)} \dots \partial x^{(s_k)}} (x)
      \right| \leq c_{85} \cdot B^{(L-1) \cdot k} \cdot A^k
      .
      \end{equation}
We have
    \begin{eqnarray*}
      \frac{\partial f_{k,i}^{(l)}}{\partial x^{(s)}}(x)
&
      =
&
      \sigma^\prime \left(\sum_{t=1}^{r} w_{k,i,t}^{(l-1)} \cdot f_{k,t}^{(l-1)}(x) + w_{k,i,0}^{(l-1)} \right)
      \cdot
      \sum_{j=1}^{r} w_{k,i,j}^{(l-1)} \cdot \frac{\partial
f_{k,j}^{(l-1)}
      }{\partial x^{(s)}}(x)
      \\
      &
      =&
      \sum_{j=1}^{r} w_{k,i,j}^{(l-1)} 
      \cdot
      \sigma^\prime \left(\sum_{t=1}^{r} w_{k,i,t}^{(l-1)} \cdot f_{k,t}^{(l-1)}(x) + w_{k,i,0}^{(l-1)} \right)
      \cdot
      \frac{\partial
f_{k,j}^{(l-1)}
      }{\partial x^{(s)}}(x)
    \end{eqnarray*}
    and
    \begin{eqnarray*}
      \frac{\partial f_{k,i}^{(1)}}{\partial x^{(s)}}(x)
      &=&
      \sigma^\prime \left(\sum_{j=1}^d w_{k,i,j}^{(0)} \cdot x^{(j)} + w_{k,i,0}^{(0)} \right)
      \cdot
       w_{k,i,s}^{(0)}.
      \end{eqnarray*}
    By the product rule of derivation we can conclude for $l>1$ that
    \begin{equation}
      \label{ple12eq2}
    \frac{\partial^k f_{k,i}^{(l)}}{\partial x^{(s_1)} \dots \partial x^{(s_k)}} (x)
    \end{equation}
    is a sum of at most $r \cdot (r+k)^{k-1}$ terms of the form
    \begin{eqnarray*}
      &&
w
\cdot
\sigma^{(s)} \left(\sum_{j=1}^{r} w_{k,i,j}^{(l-1)} \cdot f_{k,j}^{(l-1)}(x) + w_{k,i,0}^{(l-1)} \right)
\\
&&
\hspace*{3cm}
      \cdot
      \frac{\partial^{t_{1}}
f_{k,j_1}^{(l-1)}
      }{\partial x^{(r_{1,1})} \dots \partial x^{(r_{1,t_{1}})}}(x)
      \cdot
\dots
      \cdot
      \frac{\partial^{t_{s}}
f_{k,j_s}^{(l-1)}
      }{\partial x^{(r_{s,1})} \dots \partial x^{(r_{s,t_{s}})}}(x)
    \end{eqnarray*}
    where we have $s \in \{1, \dots,k\}$,
    $|w| \leq  B^{s}$ and $t_1+ \dots + t_s = k$.
    Furthermore
    \[
    \frac{\partial^k f_{k,i}^{(1)}}{\partial x^{(s_1)} \dots \partial x^{(s_k)}} (x)
    \]
    is a given by
    \[
    \prod_{j=1}^k w_{k,i,s_j}^{(0)} \cdot
      \sigma^{(k)} \left(\sum_{t=1}^d w_{k,i,t}^{(0)} \cdot x^{(t)} + w_{k,i,0}^{(0)} \right)
.
    \]
    Because of the boundedness of the derivatives of $\sigma$ we can conclude
    from (\ref{le12eq3})
    \[
    \left|
    \frac{\partial^k f_{k,i}^{(1)}}{\partial x^{(s_1)} \dots \partial x^{(s_k)}} (x)
\right|
\leq
c_{86} \cdot
A^k
    \]
    for all $k \in \N$ and $s_1, \dots, s_k \in \{1, \dots, d\}$.

    Recursively we can conclude from the above
    representation of (\ref{ple12eq2}) that we have
    \[
    \left|
    \frac{\partial^k f_{k,i}^{(l)}}{\partial x^{(s_1)} \dots \partial x^{(s_k)}} (x)
    \right|
    \leq
    c_{86,r,l,k} \cdot B^{(l-1) \cdot k} \cdot A^k.
    \]
    Setting $l=L$ we get (\ref{ple12eq1b}).

In the {\it second step} of the proof we show
\begin{equation}
	\label{ple12eq3}
	\Nu_p \left(
	\epsilon, \{T_\beta f \cdot 1_{[-\alpha,\alpha]^d} \, : \, f \in \F \}, x_1^n
	\right)
	\leq
	\Nu_p  \left(
	\frac{\epsilon}{2}, T_\beta \G\circ \Pi, x_1^n
	\right),     
\end{equation}
where $\G$ is the set of all polynomials of degree less than or equal to $k-1$
which vanish outside of $[-\alpha,\alpha]^d$
and $\Pi$ is the family of all partitions of $\mathbb{R}^d$ which consist of a partition of $[-\alpha,\alpha]^d$ into
\[
K =
\left(
\left\lceil
\frac{2 \cdot \alpha}{
  \left(
c_{87} \cdot \frac{\epsilon}{c}
  \right)^{1/k}
}
\right\rceil
\right)^d
\]
many cubes of sidelenght at most
\[
\left(c_{87}\cdot \frac{\epsilon}{c}\right)^{1/k}
\]
where $c_{87}=c_{87}(d,k)>0$ is a suitable small constant greater than zero,
and the additional set $\mathbb{R}^d\setminus [-\alpha,\alpha]^d$.

A standard bound on the remainder of a multivariate Taylor polynomial
together with (\ref{ple12eq1})
shows that for each $f_\bw$ we can find $g \in \mathcal{G}\circ \Pi$ such that
\[
|f_\bw(x)-g(x)|\leq \frac{\epsilon}{2}
\]
  holds for all $x\in [-\alpha,\alpha]^d$, which implies (\ref{ple12eq3}).

    In the {\it third step} of the proof we show the assertion of Lemma \ref{le12}.
    Since $\mathcal{G}\circ \Pi$ is a linear vector space of dimension less than or equal to \[
    c_{88} \cdot \alpha^d \cdot \left(\frac{c}{\epsilon}\right)^{d/k}
    \]
    we conclude from Theorem 9.4 and Theorem 9.5 in Gy\"orfi et al. (2002),
    	\begin{align*}
    	\mathcal{N}_p(\frac{\epsilon}{2}, T_\beta \mathcal{G} \circ \Pi, x_1^n)
    	\leq 3 \left(\frac{2e(2\beta)^p}{(\epsilon/2)^p}\log\left(\frac{3e(2 \beta)^p}{(\epsilon/2)^p}\right)\right)^{c_{88} \cdot \alpha^d \cdot \left(\frac{c}{\epsilon}\right)^{d/k}+1}.
    \end{align*}
    
\noindent    
    Together with (\ref{ple12eq3}) this implies the assertion.
    \quad  \hfill $\Box$

    \subsection{Proof of Theorem \ref{th1}}

W.l.o.g. we assume
throughout the proof that $n$ is sufficiently large and that
$\|m\|_\infty \leq \beta_n$ holds.
Let $A>0$ with $supp(X) \subseteq [-A,A]^d$.
Set
\[
\tilde{K}_n= \left\lceil c_{89} \cdot n^{\frac{d}{2p+d}} \right\rceil 
\]
and
\[
N_n= \left\lceil c_{90} \cdot n^{3+\frac{d}{2p+d}}  \right\rceil
\]
and let $\bw^*$ be a weight vector of a neural networks
where the results of $N_n \cdot \tilde{K}_n$ in parallel computed neural
networks with $L$ hidden layers and $r$ neurons per layer are
computed such that the corresponding network
\[
f_{\bw^*}(x)= \sum_{k=1}^{N_n \cdot \tilde{K}_n} (\bw^*)_{1,1,k} \cdot f_{\bw^*,k,1}^{(L)}(x)
\]
satisfies
\begin{equation}
\label{pth1eq3}
\sup_{x \in [-A,A]^d} |f_{\bw^*}(x)-m(x)| \leq \frac{c_{91}}{\tilde{K}_n^{p/d} }
\end{equation}
and
\[
|(\bw^*)_{1,1,k}| \leq \frac{c_{92} \cdot  \tilde{K}_n^{(q+d)/d}}{N_n}
\quad (k=1, \dots, N_n \cdot \tilde{K}_n).
\]
Note that such a network exists according Theorem \ref{se4th2}
if we repeat in the outer sum of the function space $\HH$
each of the $f_k$'s in Theorem \ref{se4th2} $N_n$--times
with outer weights divided by $N_n$.
Set
\[
\epsilon_n = \frac{c_{93}}{ n \cdot \sqrt{N_n \cdot \tilde{K}_n}}
\geq
\frac{c_{94}}{n^{4}}
\]
where the last inequality holds because of $p \geq 1/2$.
   Let $A_n$ be the event that firstly the weight vector $\bw^{(0)}$
            satisfies
            \[
            | (\bw^{(0)})_{j_s,k,i}^{(l)}-(\bw^*)_{s,k,i}^{(l)}| \leq \epsilon_n
            \quad \mbox{for all } l \in \{0, \dots, L-1\},
            s \in \{1, \dots, N_n \cdot \tilde{K}_n \}
            \]
            for some pairwise distinct $j_1, \dots, j_{N_n \cdot \tilde{K}_n}
            \in \{1, \dots, K_n\}$
and such that secondly
\[
\max_{i=1, \dots, n} |Y_i| \leq \sqrt{\beta_n}
\]
holds.

We
decompose the  $L_2$ error of $m_n$ in a sum of several terms.
Set
\[
m_{\beta_n}(x)=\EXP\{ T_{\beta_n} Y | X=x \}.
\]
We have
\begin{eqnarray*}
&&
\int | m_n(x)-m(x)|^2 \PROB_X (dx)
\\
&&
=
\left(
\EXP \left\{ |m_n(X)-Y|^2 | \D_n \right\}
-
\EXP \{ |m(X)-Y|^2\}
\right)
\cdot 1_{A_n}
+
\int | m_n(x)-m(x)|^2 \PROB_X (dx)
\cdot 1_{A_n^c}
\\
&&
=
\Big[
\EXP \left\{ |m_n(X)-Y|^2 | \D_n \right\}
-
\EXP \{ |m(X)-Y|^2\}
\\
&&
\hspace*{2cm}
- \left(
\EXP \left\{ |m_n(X)-T_{\beta_n} Y|^2 | \D_n \right\}
-
\EXP \{ |m_{\beta_n}(X)- T_{\beta_n} Y|^2\}
\right)
\Big] \cdot 1_{A_n}
\\
&&
\quad +
\Big[
\EXP \left\{ |m_n(X)-T_{\beta_n} Y|^2| \D_n \right\}
-
\EXP \{ |m_{\beta_n}(X)- T_{\beta_n} Y|^2\}
\\
&&
\hspace*{2cm}
-
2 \cdot \frac{1}{n} \sum_{i=1}^n
\left(
|m_n(X_i)-T_{\beta_n} Y_i|^2
-
|m_{\beta_n}(X_i)- T_{\beta_n} Y_i|^2
\right)
\Big] \cdot 1_{A_n}
\\
&&
\quad
+\Big[
2 \cdot \frac{1}{n} \sum_{i=1}^n
|m_n(X_i)-T_{\beta_n} Y_i|^2
-
2 \cdot \frac{1}{n} \sum_{i=1}^n
|m_{\beta_n}(X_i)- T_{\beta_n} Y_i|^2
\\
&&
\hspace*{2cm}
- \left(
2 \cdot \frac{1}{n} \sum_{i=1}^n
|m_n(X_i)-Y_i|^2
-
2 \cdot \frac{1}{n} \sum_{i=1}^n
|m(X_i)- Y_i|^2
\right)
\Big] \cdot 1_{A_n}
\\
&&
\quad
+
\Big[
2 \cdot \frac{1}{n} \sum_{i=1}^n
|m_n(X_i)-Y_i|^2
-
2 \cdot \frac{1}{n} \sum_{i=1}^n
|m(X_i)- Y_i|^2
\Big] \cdot 1_{A_n}
\\
&&
\quad
+
\int | m_n(x)-m(x)|^2 \PROB_X (dx)
\cdot 1_{A_n^c}
\\
&&
=: \sum_{j=1}^5 T_{j,n}.
\end{eqnarray*}
In the remainder of the proof we bound
\[
\EXP T_{j,n}
\]
for $j \in \{1, \dots, 5\}$.

In the {\it first step of the proof} we show
\[
\EXP T_{j,n} \leq c_{95} \cdot \frac{\log n}{n} \quad
\mbox{for } j \in \{1,3\}.
\]
This follows as in the proof of Lemma 1 in Bauer and Kohler (2019).

In the {\it second step of the proof} we show
\[
\EXP T_{5,n} \leq c_{96} \cdot \frac{(\log n)^2}{n}.
\]

The definition of $m_n$ implies $\int |m_n(x)-m(x)|^2 \PROB_X (dx) \leq
4 \cdot c_3^2 \cdot (\log
n)^2$, hence it suffices to show
\begin{equation}
\label{pth1eq1}
\PROB(A_n^c) \leq \frac{c_{97}}{n^2}.
\end{equation}
To do this, we consider a sequential choice of the weights of the
$K_n$ fully connected neural networks. The probability that the weights in
the first of these networks differ in all components at most by $\epsilon_n$
from $(\bw^*)_{1,i,j}^{(l)}$ $(l=0, \dots, L-1)$ is
for large $n$ bounded from below by
\begin{eqnarray*}
  &&
  \left( \frac{c_{94}}{2 \cdot c_{1}  \cdot n^{4}}
\right)^{r \cdot (r+1) \cdot (L-1)}
\cdot
\left(
\frac{1}{2 \cdot c_2 \cdot (\log n) \cdot n^\tau \cdot n^{4}}
\right)^{r \cdot (d+1)}
\\
&&
\geq
n^{-r \cdot (r+1) \cdot (L-1) \cdot 4- r \cdot 4 \cdot (d+1) - r \cdot 4  \cdot \tau - 0.5}.
\end{eqnarray*}
Hence probability that none of the first $n^{
r \cdot (r+1) \cdot (L-1) \cdot 4+ r \cdot 4 \cdot (d+1)+ r \cdot 4 \cdot \tau +1}$ neural networks satisfies this condition is for large $n$
 bounded above by
 \begin{eqnarray*}
   &&
(1 -  n^{-r \cdot (r+1) \cdot (L-1) \cdot 4- r \cdot 4 \cdot (d+1)- r \cdot 4 \cdot \tau -0.5}) ^{n^{r \cdot (r+1) \cdot (L-1) \cdot 4 + r \cdot 4 \cdot (d+1) + r \cdot 4 \cdot \tau +1}}\\
&&\leq
\left(\exp \left(
-  n^{-r \cdot (r+1) \cdot (L-1) \cdot 4- r \cdot 4 \cdot (d+1) - r \cdot 4 \cdot \tau -0.5}
\right)
\right) ^{n^{r \cdot (r+1) \cdot (L-1) \cdot 4+ r \cdot 4 \cdot (d+1) + r \cdot 4 \cdot \tau +1}}
\\
&&=
\exp( -  n^{0.5}).
\end{eqnarray*}
 Since we have $K_n \geq n^{r \cdot (r+1) \cdot (L-1) \cdot 4+ r \cdot 4 \cdot (d+1) +
   r \cdot 4 \cdot \tau +1} \cdot N_n \cdot \tilde{K}_n$
 for $n$ large we can successively
use the same construction for all of $N_n \cdot \tilde{K}_n$ weights and we can conclude:
The probability that there exists $k \in \{1, \dots, N_n \cdot \tilde{K}_n\}$
such that
none of the $K_n$ weight vectors of the fully
connected neural network differs by at most $\epsilon_n$ from
$((\bw^*)_{i,j,k}^{(l)})_{i,j,l}$ is for large $n$ bounded from above by
\begin{eqnarray*}
&&
N_n \cdot \tilde{K}_n \cdot \exp( -  n^{0.5})
\leq  c_{98} \cdot n^{5} \cdot \exp( - n^{0.5})
\leq \frac{c_{99}}{n^2}.
\end{eqnarray*}
This implies for large $n$
\begin{eqnarray*}
\PROB(A_n^c)
&\leq&
\frac{c_{99}}{n^2}
+
\PROB\{ \max_{i=1, \dots, n} |Y_i| > \sqrt{\beta_n}
\}
 \leq
\frac{c_{99}}{n^2}
+
n \cdot\PROB\{  |Y| > \sqrt{\beta_n}
\}
\\
& \leq &
\frac{c_{99}}{n^2}
+
n \cdot
\frac{\EXP\{ \exp(c_4 \cdot Y^2)}{ \exp( c_4 \cdot \beta_n)}
\leq
\frac{c_{97}}{n^2},
\end{eqnarray*}
where the last inequality holds because of (\ref{th1eq1})  and
$c_3 \cdot c_4 \geq 2$.

Let $\epsilon >0$ be arbitrary.
In the {\it third step of the proof} we show
\[
\EXP T_{2,n} \leq
c_{100} \cdot
\frac{ n^{\tau \cdot d + \epsilon}}{n}
.
\]
Let $\W_n$ be the set of all weight vectors
$(w_{i,j,k}^{(l)})_{i,j,k,l}$ which satisfy
\[
| w_{1,1,k}^{(L)}| \leq c_{101} \quad (k=1, \dots, K_n),
\]
\[
|w_{i,j,k}^{(l)}| \leq c_{102} \quad (l=1, \dots, L-1)
\]
and
\[
|w_{i,j,k}^{(0)}| \leq (c_2+c_{103}) \cdot (\log n) \cdot n^\tau.
\]
By  Lemma \ref{le4}, Lemma \ref{le5} and Lemma  \ref{le6} we can
conclude that on $A_n$ we have
\begin{equation}
\label{pth1eq2}
\| \bw^{(t)}-\bw^{(0)}\| \leq c_{104} \quad (t=1,
\dots, t_n).
\end{equation}
This follows from the fact that on $A_n$ we have
\[
F_n(\bw^{(0)})
=
\frac{1}{n} \sum_{i=1}^n Y_i^2 \leq \beta_n
\]
and that
\[
2 \cdot t_n \cdot \lambda_n \cdot \beta_n
\leq
c_{105}.
\]
Together with the initial choice of $\bw^{(0)}$ this implies that on
$A_n$
we have
\[
\bw^{(t)} \in \W_n \quad (t=0, \dots, t_n).
\]
Hence, for any $u>0$ we get
\begin{eqnarray*}
&&
\PROB \{ T_{2,n} > u \}
\\
&&
\leq
\PROB \Bigg\{
\exists f \in \F_n :
\EXP \left(
\left|
\frac{f(X)}{\beta_n} - \frac{T_{\beta_n}Y}{\beta_n}
\right|^2
\right)
-
\EXP \left(
\left|
\frac{m_{\beta_n}(X)}{\beta_n} - \frac{T_{\beta_n}Y}{\beta_n}
\right|^2
\right)
\\
&&\hspace*{3cm}-
\frac{1}{n} \sum_{i=1}^n
\left(
\left|
\frac{f(X_i)}{\beta_n} - \frac{T_{\beta_n}Y_i}{\beta_n}
\right|^2
-
\left|
\frac{m_{\beta_n}(X_i)}{\beta_n} - \frac{T_{\beta_n}Y_i}{\beta_n}
\right|^2
\right)
\Bigg\}
\\
&&\hspace*{2cm}
> \frac{1}{2} \cdot
\left(
\frac{u}{\beta_n^2}
+
\EXP \left(
\left|
\frac{f(X)}{\beta_n} - \frac{T_{\beta_n}Y}{\beta_n}
\right|^2
\right)
-
\EXP \left(
\left|
\frac{m_{\beta_n}(X)}{\beta_n} - \frac{T_{\beta_n}Y}{\beta_n}
\right|^2
\right)
\right),
\end{eqnarray*}
where
\[
\F_n = \left\{ T_{\beta_n} f_\bw \quad : \quad \bw \in \W_n \right\}.
\]
By Lemma \ref{le12} we get
\begin{eqnarray*}
&&
\Nu_1 \left(
\delta , \left\{
\frac{1}{\beta_n} \cdot f : f \in \F_n
\right\}
, x_1^n
\right)
\leq
\Nu_1 \left(
\delta \cdot \beta_n , \F_n
, x_1^n
\right)
\\
&&
\leq
\left(
\frac{ c_{106}}{\delta}
\right)^{
c_{107} \cdot (\log n)^{d}  n^{\tau \cdot d} \cdot (c_{108})^{(L-1) \cdot d}
\cdot
   \left(\frac{K_n \cdot c_{109}}{\beta_n \cdot \delta}\right)^{d/k} + c_{110}
  }.
\end{eqnarray*}
By choosing $k$ large enough we get for $\delta>1/n^2$
\[
\Nu_1 \left(
\delta , \left\{
\frac{1}{\beta_n} \cdot f : f \in \F_n
\right\}
, x_1^n
\right)
\leq
c_{111} \cdot n^{ c_{112} \cdot n^{\tau \cdot d + \epsilon/2}}.
\]
This together with Theorem 11.4 in Gy\"orfi et al. (2002) leads for $u
\geq 1/n$ to
\[
\PROB\{T_{2,n}>u\}
\leq
14 \cdot
c_{111} \cdot n^{ c_{112} \cdot n^{\tau \cdot d +  \epsilon/2}}
\cdot
\exp \left(
- \frac{n}{5136 \cdot \beta_n^2} \cdot u
\right).
\]
For $\epsilon_n \geq 1/n$ we can conclude
\begin{eqnarray*}
\EXP \{ T_{2,n} \}
& \leq &
\epsilon_n + \int_{\epsilon_n}^\infty \PROB\{ T_{2,n}>u \} \, du
\\
& \leq &
\epsilon_n
+
14 \cdot
c_{111} \cdot n^{ c_{112} \cdot n^{\tau \cdot d +
    \epsilon/2}}
\cdot
\exp \left(
- \frac{n}{5136 \cdot \beta_n^2} \cdot \epsilon_n
\right)
\cdot
\frac{5136 \cdot \beta_n^2}{n}.
\end{eqnarray*}
Setting
\[
\epsilon_n = \frac{5136 \cdot \beta_n^2}{n}
\cdot
c_{112}
\cdot
 n^{\tau \cdot d +
    \epsilon/2}
\cdot \log n
\]
yields the assertion of the fourth step of the proof.

In the {\it fourth step of the proof} we show
\begin{eqnarray*}
&&
  \EXP \{ T_{4,n} \} 
  \leq c_{113} \cdot n^{- \frac{2p}{2p+d}}.
  \end{eqnarray*}
Using
\[
|T_{\beta_n} z - y| \leq |z-y|
\quad \mbox{for } |y| \leq \beta_n
\]
we get
\begin{eqnarray*}
&&
 T_{4,n}/2
\\
&&
=
\Big[ \frac{1}{n} \sum_{i=1}^n
|m_n(X_i)-Y_i|^2
-
 \frac{1}{n} \sum_{i=1}^n
|m(X_i)- Y_i|^2
\Big] \cdot 1_{A_n}
\\
&&
\leq
\Big[
\frac{1}{n} \sum_{i=1}^n
|f_{\bw^{(t_n)}}(X_i)-Y_i|^2
-
 \frac{1}{n} \sum_{i=1}^n
|m(X_i)- Y_i|^2
\Big] \cdot 1_{A_n}
\\
&&
\leq
\big[ F_n(\bw^{(t_n)})
-
 \frac{1}{n} \sum_{i=1}^n
|m(X_i)- Y_i|^2
\Big] \cdot 1_{A_n}.
\end{eqnarray*}
On $A_n$ we have 
\begin{eqnarray*}
\| \bw^{*} - \bw^{(0)} \|^2
&\leq&
\sum_{k=1}^{K_n} |(\bw^{*})_{1,1,k}^{(L)}|^2 + N_n \cdot \tilde{K}_n
   \cdot L \cdot (r \cdot
   (r+d))^L \cdot \epsilon_n^2
\\
&\leq&
N_n \cdot \tilde{K}_n \cdot \left(
\frac{c_{114} \cdot \tilde{K}_n^{(q+d)/d}}{N_n}
\right)^2
+
c_{115} \cdot N_n \cdot \tilde{K}_n
   \cdot \epsilon_n^2 \leq \frac{c_{116}}{n^2}. 
\end{eqnarray*}

Application of Theorem \ref{se4th1}
yields
\begin{eqnarray*}
  &&
  T_{4,n}/2
  \\
  &&
 \leq
 \Bigg(
\frac{1}{n} \sum_{i=1}^n
|f_{\bw^{*}}(X_i)-Y_i|^2
+ c_{117} \cdot (\log n) \cdot n \cdot \|\bw^*-\bw^{(0)}\|^2
+ c_{118} \cdot \frac{\log n}{n}
\\
&&
\quad \quad
-
 \frac{1}{n} \sum_{i=1}^n
|m(X_i)- Y_i|^2
 \Bigg)
 \cdot 1_{A_n}
 \\
 &&
 \leq
 \Bigg(
\frac{1}{n} \sum_{i=1}^n
|f_{\bw^{*}}(X_i)-Y_i|^2
+ c_{119} \cdot (\log n) \cdot n \cdot \frac{1}{n^2}
+ c_{120} \cdot \frac{\log n}{n}
\\
&&
\quad \quad
-
 \frac{1}{n} \sum_{i=1}^n
|m(X_i)- Y_i|^2
 \Bigg)
 +
 \frac{1}{n} \sum_{i=1}^n
|m(X_i)- Y_i|^2
 \cdot 1_{A_n^c}. 
\end{eqnarray*}

            Hence
            \begin{eqnarray*}
              &&
              \EXP\{T_{4,n}\}
                \\
                &&
                \leq
                2 \cdot \int |f_{\bw^*}(x)-m(x)|^2 \PROB_X(dx)
                + c_{119} \cdot (\log n) \cdot n \cdot \frac{1}{n^2}
+ c_{120} \cdot \frac{\log n}{n}
\\
&&
\quad
+
\sqrt{
  \EXP \left\{
  \left|
 \frac{1}{n} \sum_{i=1}^n
|m(X_i)- Y_i|^2
  \right|^2
  \right\}
}
\cdot \sqrt{\PROB(A_n^c)}
\\
&&
\leq
                2 \cdot \int |f_{\bw^*}(x)-m(x)|^2 \PROB_X(dx)
                + c_{121} \cdot \log n \cdot \frac{1}{n}.
              \end{eqnarray*}
Application of (\ref{pth1eq3}) yields the assertion.
\hfill $\Box$



\end{document}